\documentclass[english,fleqn]{article}
\usepackage[left=1.5cm,right=1.5cm,top=1.5cm,bottom=1.5cm]{geometry}
\usepackage[T1]{fontenc}
\usepackage[utf8]{inputenc}
\usepackage{geometry}
\usepackage{authblk}
\usepackage[cmex10]{amsmath}
\usepackage{amssymb}
\usepackage{graphicx}
\usepackage{babel}
\usepackage{setspace}
\usepackage{supertabular}
\usepackage{multirow,booktabs}
\usepackage{pgfplotstable}
\pgfplotsset{compat=1.14}
\usepackage{lineno}
\usepackage[numbers]{natbib}
\numberwithin{equation}{section}
\usepackage{url}
\usepackage{breakurl}
\usepackage[breaklinks]{hyperref}
\providecommand{\keywords}[1]{\textbf{\textit{Index terms---}} #1}
\usepackage[section]{algorithm}
\usepackage{algorithmic}
\usepackage{color}
\usepackage[color=green!40]{todonotes}
\def\comment[#1]#2{\todo[inline,author=#1]{#2}}
\usepackage[authormarkup=none]{changes}

\begin{document}

\title{A hybrid Benders decomposition and bees algorithm matheuristic approach
to transmission expansion planning considering energy storage}

%

\author[1]{C.A.G. MacRae \thanks{Corresponding author; Email: \texttt{cmacrae@nust.na}}}
\author[2]{M. Ozlen}
\author[3]{A.T. Ernst}

\affil[1]{\small{Department of Computer Science, Namibia University of Science
and Technology }}
\affil[2]{\small{School of Science, RMIT University}}
\affil[3]{\small{School of Mathematical Sciences, Monash University}}

\maketitle

\begin{abstract}

  This paper introduces a novel hybrid optimisation algorithm that combines
  elements of both metaheuristic search and integer programming. This new
  matheuristic combines elements of Benders decomposition and the Bees
  Algorithm, to create the Bee-Benders Hybrid Algorithm (BBHA) which retains
  many of the advantages both of the methods. Specifically it is designed to be
  easily parallelizable, to produce good solutions quickly while still retaining
  a guarantee of optimality when run for a sufficiently long time. The algorithm
  is tested using a transmission network expansion and energy storage planning
  model, a challenging and very large scale mixed integer linear programming
  problem. Transmission network planning problems are already difficult on their
  own. When including the planning for storage systems in the network, the
  variation of demand over time has to be taken into account significantly
  increasing the size and difficulty of the optimization problem.  The BBHA is
  shown to be highly effective hybrid matheuristic algorithm that performs at
  least as well as either Benders decomposition or the Bees Algorithm where
  these are effective on their own, and significantly improves upon the
  individual approaches where neither component part has a pronounced advantage.
  While the paper demonstrates the effectiveness in terms of the concrete
  electricity network planning problem, the algorithm could be readily applied
  to any mixed integer linear program, and is expected to work particularly well
  whenever this has a structure that is amenable to Benders decomposition.   

\end{abstract}

\keywords{hybrid heuristic, optimization, power transmission, energy storage}

\section{Introduction}
\label{sec:intro}

The need to solve large scale mixed integer programming problems arises in many 
applications. In this paper the planning of electricity networks to cope with
renewable generation has been used to both motivate the need for the proposed
new algorithm and to evaluate its effectiveness.

Integrating renewable energy generation, especially variable generators such as
wind and solar, into the electrical transmission network is a considerable
design challenge currently facing network planners. For example, a recent
blackout in South Australia saw 315MW of wind generation disconnect from the
grid amid voltage dips and loss of load \citep{chang_south_2016}, recently a 100
mega-watt battery was installed by Tesla to prevent this from happening again
\citep{keck2019impact}. Correspondingly, there has been a renewed interest in
electricity network planning problems \citep{gu_coordinating_2012}. One such
problem is the transmission expansion planning problem (TEP). In this problem
the objective of the planning is minimize the investment and operational costs
of the network while meeting a set of operational constraints, for example,
generation, demand, geographical, and environmental constraints
\citep{latorre_classification_2003}. Transmission lines require a huge initial
investment but have long life whereas storage systems have a short life but can
be installed very quickly and be expanded gradually, \citep{Qiu2017}.

Hydro is the most common type of transmission scale Energy Storage System (ESS),
but its feasibility is determined by climate, geography, and environmental
constraints. Batteries have also been successfully deployed to smooth the 5
minute ramp rate of a wind farm \citep{wood_integrating_2012} and given their
high power and energy capacities compressed air technologies remain viable but
expensive \citep{zhao_review_2015}. The transmission network expansion and
energy storage planning (TESP) model considers a generic ESS, the primary
purposes of which are transmission upgrade deferral and demand shifting.
Transmission upgrade deferral occurs when need for additional or larger capacity
transmission lines is avoided, in this case by using ESS located near sites of
generation or demand, to store energy and release it at a steady rate over time.
Subtly different, demand shifting entails using stored energy generated in a
prior time interval to meet demand in the current time interval. Storage
facilities are also considered as decision makers dealing with risks arising
from the costs, environmental impact and supply issues \citep{Sheikh2015}.

The TEP and related problems are often modeled as a mixed integer nonlinear
program (MINLP), or in a correspondent disjunctive mixed integer linear
programming (MIP) form. An overview of the standard models and test systems is
given in \citep{romero_test_2002}.

Advances in commercial solver technology mean that simpler linear models of
small networks can generally be solved to optimality within a few minutes.
However, a considerable body of research is dedicated to solving larger or more
complicated instances. Novel approaches to these problems include branch and
bound with a GRASP meta-heuristic \citep{bahiense_mixed_2001},
Projection-Adapted Cross Entropy \citep{eshragh_projection-adapted_2011}, and
particle swarm optimization \citep{aghaei_distribution_2014}. Often TEP problems
can be decomposed into investment and operational subproblems. Benders
decomposition with alternately continuous or discrete decision variables in the
master (investment) problem, and DC approximation or transportation operational
subproblems is investigated in \citep{pereira_decomposition_1985}. Fencing
constraints and additional constraints on new paths are shown to substantially
reduce the number of iterations when added to the master problem
\citep{haffner_branch_2000}, and adding Gomory cuts evaluated from the master
problem to a Benders decomposition of a linear disjunctive MIP model is shown to
result in significant CPU time savings \citep{binato_new_2001}. More recently,
local branching is used to accelerate the Benders decomposition of a TEP problem
in \citep{dilwali_transmission_2016}. A useful survey of the literature is given
in \citep{sorokin_algorithms_2012}.

Where the modeling becomes more complex and computationally demanding,
meta-heuristic approaches have been shown to produce good results. If the
transmission expansion planner is concerned only with determining a final
network plan, the planning is considered static, whereas dynamic planning
involves the determination of one or more intermediate plans over multiple
periods. A specialized genetic algorithm is shown to produce good solutions for
coordinated, multistage planning problems\citep{escobar_multistage_2004}. A
Differential evolution algorithm is used to solve a similar problem in
\citep{sum-im_differential_2009}.

Another common complication to the planning is (n-1) redundancy. Systems
operating under this scheme must not shed load if a single component, in this
context a circuit, fails. An adaption of the Chu-Beasley\citep{chu_genetic_1997}
genetic algorithm was used to solve a TEP model with (n-1) security constraints
in \citep{de_j_silva_transmission_2005}.


Some recent works have considered incorporating energy storage in transmission
networks, however the time dimension is largely ignored allowing the ESS to
behave as an alternative type of generation \citep{hu_transmission_2012}. A pair
of linear programming models that take into account both variable and
dispatchable generation, as well as energy storages are compared in
\citep{clack_linear_2015}.

Locating and sizing small scale ESS in distribution networks has been approached
using a modified particle swarm optimization (PSO) to optimize a multi-period
design problem\citep{sedghi_distribution_2013}, and using a genetic algorithm
combined with simulated annealing to plan a low voltage network with high solar
photovoltaic generation \citep{crossland_planning_2014}.

For solving large scale optimisation problems such as encountered in network
planning, there are two categories of methods that are commonly used: (1)
meta-heuristics that aim to get good solutions quickly but have no guarantee of
optimality, and (2) exact methods such as integer programming and constraint
programming that at least in principle produce optimal solution but may take too
long to run in practice. Over the last decade there has been an increasing
interest in hybrid methods that combined elements taken from both worlds. While
a complete review of such methods is beyond the scope of this paper, a few
examples show the variety of approaches tried in this category.
\citep{puchinger_metaboosting2010} developed a method they called MetaBoosting
to use metaheuristics to improve column generation or cut separation in an
integer programming framework. Using Ant Colony Optimisation with Constraint
Propagation has been shown to be effective in
\cite{thiruvady_constraint-based_2013}. Lagrangian relaxation has been combined
with Particle Swarm
Optimisation~\cite{ernst_hybrid_2010,gomez-iglesias_scalable_2013} or Ant Colony
Optimisation~\cite{thiruvady_lagrangian-aco_2014}. Genetic programming has been
combined with MIP though largely as a high-level tuning mechanism for the large
number of parameters available in a MIP solver~\cite{kostikas_genetic_2004}.
The rise of these types of methods that, in principle, can operate on any
integer linear programming problem, has led to the term \emph{matheuristics}
being used to describe such methods that combine mathematical programming with
heuristics \cite{boschetti_matheuristics:_2009}.

Of particular interest from the point of view of this paper are methods based on
combining Benders decomposition with metaheuristic algorithms.  Benders
decomposition breaks large MIPs into a master problem and one or more
subproblems. The subproblems are used to evaluate and test the feasibility of
solutions proposed by an optimisation process solving the master problem.
Information from the subproblems is also fed back to the master in terms of
additional constraints (cuts) generated from the subproblem solutions.
\citep{poojari_improving_2009} created a method based on Benders decomposition
where the master problem is always solved with a Genetic Algorithm (GA) while
the subproblem is solved using Linear Programming (LP) as in the standard
Benders approach. When tested on general MIPs this proved to be more effective
than using the Benders method on its own but often less effective than simply
using the CPLEX solver without any decomposition. A little earlier
\cite{sirikum_new_2007} independently developed a nearly identical GA-Benders
hybrid and applied it to a problem in power generation expansion planning, a
similar problem to the one considered here but only looking at power generation
investment without any consideration of network transmission or energy storage.
  
In this paper, we present a hybrid exact/meta-heuristic algorithm that melds
Benders decomposition and a Bees Algorithm (BA)\citep{pham_bees_2009} inspired
approach. Unlike previous proposed matheuristics, this method retains the
ability to generate provably optimal solutions. The ideas of multiple bees that
have different functions, scouts and workers, is used to balance diversification
and intensification in the solution of the master problem.  Multiple parallel
optimisation processes are used to speed up the search, that learn from each
other not just through the exchange of good feasible solutions but also by
sharing dual (cut) information obtained when solving Benders subproblems. Using
the transmission network expansion and energy storage planning model (TESP) to
test the model, we show the Bee-Benders hybrid algorithm (BBHA) to be an
effective hybrid algorithm that exhibits equivalent performance to its component
parts in the segments of the problem domain where those parts are strongest, and
significantly improves upon the individual approaches where neither component
part has a pronounced advantage.

The rest of this paper is organized as follows. The Bee-Benders algorithm is
introduced in Section~\ref{sec:BBHA}. A MIP formulation of the TEP
with storage model is given in Section~\ref{sec:model}. Numerical results,
in which the algorithm is evaluated using the Brazilian 46-bus and Colombian
93-bus test systems are discussed in Section~\ref{sec:num_res}. We conclude in
Section~\ref{sec:conclusion}.

\section{The Bee-Benders Hybrid Algorithm}
\label{sec:BBHA}

\subsection{The Bees Algorithm}
Here we present a hybrid exact/meta-heuristic algorithm that combines Benders
decomposition with an approach inspired by the Bees Algorithm.
There are many variants of optimisation metaheuristics inspired by the
behaviour of bees (see for example \cite{karaboga2014} for a review of one of
the alternatives, the Artificial Bee Colony optimisation). Here we will follow
the Bees Algorithm as proposed by \citep{pham_bees_2009} and
\citep{pham2007using}.

In the most basic form, the algorithm comprises two phases: global search, and
local search. A pseudo-code description of this metaheuristic has been provided
in Algorithm~\ref{alg:BA}. Each solution in the solution space is referred to as
a flower in the terminology of the Bees Algorithm, and the local neighbourhood
of a solution is called a flower patch. In the initialisation phase, ``scout''
bees leave the hive and fly to a random flower. The fitness of the flower is
evaluated and the scout bees return to the hive. During the local search phase,
the scouts who discover the $ne$ elite and the $nb$ best flowers (solutions)
recruit ``worker'' bees to explore their respective flower patches, that is, the
flowers in the neighbourhood of those the scouts discovered.  Recruited worker
bees fly to a random flower within the flower patch and evaluate its fitness.
The fittest flower from the elite and best flower patches are combined with the
fittest new flowers discovered scouts during the subsequent global search phase
to produce a new pool of elite and best solutions for further local exploration.
Stopping conditions may include time, the number of iteration, or a test for
convergence.

\begin{algorithm}[htb!]
  \label{alg:BA}
  \caption{Bees Algorithm}
  \begin{algorithmic}[1]{} 
    \REQUIRE{ $ns:$ no. of scout bees}
    \REQUIRE{ $nre\geq nrb:$ no. of recruited bees per elite / best sites}
    \REQUIRE{ $ne \leq nb:$ no. of elite / best sites, $nb\leq ns$}
    \REQUIRE{ $ngh, stlim:$ initial size of neighbourhood \& stagnation limit}
    \FOR{$b=1,\ldots,ns$}{
      \STATE Generate a random initial solution to include in the set of sites $S$
    }\ENDFOR
    \WHILE{ not out of time}{
      \STATE Evaluate the fitness (objective) of all sites in $S$
      \STATE Let $S=E\cup B\cup R$ where $|E|=ne$ and $|B|=nb$ with $R$ having
      lower fitness than the others.
      \FORALL{solutions (sites) $s \in E\cup B$}{
        \STATE Evaluate $ne$ (or $nb$) solutions in the neighbourhood of $s$
        if $s\in E$ (resp. $s\in B$)
        \IF{ better solution found}{
          \STATE Replace $s$ with the best solution found.
        }\ELSE{
          \STATE Reduce the neighbourhood size $ngh$
          \STATE delete site $s$ if $stlim$ iterations without improvement
        }\ENDIF
      }\ENDFOR
      \STATE Let $S:=E\cup B$ and add random solutions until $|S|=ns$.
    }\ENDWHILE
  \end{algorithmic}

\end{algorithm}

This can be thought as a multi-start local search algorithm which always works
on a subset of best known solutions, with more effort expended on the local
search in the neighbourhood of the elite solutions than the remaining solutions.
It should be noted that this algorithm can be parallelized in a fairly straight
forward manner by carrying out the search in each flower patch (neighbourhood of
an elite or best solution) in parallel. Also the ``scouting'', that is,
generation of new random solutions, can be carried out independently.

The BA has been applied to numerous combinatorial optimization problems such as
the generalized assignment problem\citep{ozbakir_bees_2010} and machine
scheduling \citep{pham2007using}, and as also shown value for applied industrial
applications such as crack detection of beam-type structures
\citep{moradi_application_2011}.

We have chosen to combine the BA with Benders decomposition because it has been
shown to perform at least as well as standard evolutionary approaches, to be
less sensitive to tuning parameters than other swarm approaches such as PSO, and
yet retains an extraordinary simplicity of implementation \citep{pham_bees_2009}.

\subsection{Benders Decomposition}

Benders decomposition is a technique that allows a large, intractable problem,
such as the TESP model described in Section~\ref{sec:model}, to be divided into
more tractable component parts \citep{benders_partitioning_1962}.  The first
part is called the \emph{master} problem and consists of a MIP that includes all
of the integer variables and any applicable continuous variables. The second
part consists of one or more \emph{subproblems}, that collectively contain the
remaining continuous variables \citep{geoffrion_generalized_1972}. Typically
both master and subproblem(s) have not only significantly fewer variables than
the original problem but also far fewer constraints, making these much easier to
solve using a MIP solver.  The master problem is solved to yield a candidate
solution which is used to fix the complicating variables that would otherwise be
present in the subproblem.  The dual solution of the subproblem is used to
produce a feasibility or optimality cut to be added to the master, and the
master problem is solved again. This iterative procedure continues until it no
further cuts are necessary. Note that in each iteration a candidate solution is
evaluated to determine if it is both feasible and better than any seen
previously. The addition of the cuts ensures that the master problem does not
repeatedly generate the same candidate solution. For problems where only
optimality cuts are possible, as is the case in this paper, the collection of
cuts represent a piecewise linear approximation to the objective function. This
approximation is iteratively refined and improved as more cuts are added.

Benders decomposition has been applied to numerous optimization problems such as
the fixed charge network design problem\citep{costa_survey_2005}, the unit
commitment problem\citep{ma_transmission_1997}, the network-constrained unit
commitment problem \citep{wu_accelerating_2010}, and the scheduling of crude oil
in an oil refinery \citep{saharidis_accelerating_2010}. It is also proven quite
effective on multi-stage stochastic energy planning problems
\citep{Pereira1991,Rebennack2016} and CCHP-microgrid operation involving battery
storage \citep{Marino2018}.

\subsection{The hybrid method}
The Bee-Benders hybrid algorithm (BBHA) is a hybrid of Benders decomposition and
a local search phase that is largely based on the Bees Algorithm. The algorithm
operates on large MIP that has been decomposed in a manner suitable for Benders
decomposition. The master problem contains binary variables
representing certain investment decisions, and an LP subproblem containing
continuous variables and largely operational constraints. The particulars of the
mathematical model detailed in Section~\ref{sec:model}.

\subsection{The BBHA in detail}
\label{subsec:details}

As with the BA, the algorithm comprises global search and local search phases.
The global search phase commences as a conventional Benders decomposition using
a ``single tree'' master approach. Lazy constraint callbacks are used to
separate Benders cuts, as opposed to solving the master problem to optimality at
each iteration. This single tree branch and bound search fulfills the role of
the ``scout'' bee in the BA algorithm, ensuring that eventually the whole
solution space is searched and the method can never remain stuck in a local
optimum. Meanwhile, an initial set of random solutions is generated for
exploration during the local search phase.

During the local search phase, ``worker'' bees (henceforth known as ``workers'')
explore the local neighbourhood (subsequently referred to as a ``site'') of each
solution, by estimating the fitness a subset of solutions using a
matheuristic based on the set of known Benders cuts. The most heuristically
promising solution discovered at the site is selected for full evaluation of the
LP.

As with the BA, the fittest solution from both the elite and best sites are
combined with the incumbent solution of the Benders decomposition to produce a
new pool of elite and best solutions for further local search. The algorithm
iterates in this way until stopping condition is met, or the Benders
decomposition finds and proves the optimal solution. A pseudo-code description
of the BBHA is given in Algorithm~\ref{alg:BBHA}. It should be noted that while
the algorithm description is essentially serial, the intention is for each
``bee'' to be executing as a parallel process that is carrying out its search
independent of the other processes. The BBHA is discussed in greater detail
below.


\begin{algorithm}[htb!]
  \caption{Bee-Benders hybrid algorithm}
  \begin{algorithmic}[1]{} 
  \label{alg:BBHA}
    \REQUIRE{ $ne \leq nb:$ no. of elite / best sites}
    \REQUIRE{ $nre\geq nrb:$ no. of recruited workers per elite / best sites}
    \REQUIRE{ $ngh > 0 :$ maximum Hamming distance comprising a neighbourhood}
    \FOR{$b=1,\ldots,nre+nrb$}{
      \STATE Generate a random initial solution to include in the set of sites $S$
    }\ENDFOR
    \STATE Let $C$ be the set of Benders cuts
    \STATE Begin solving Benders decomposition (MIP), separating each Benders cut $c$ such that $C := C \cup \{c\}$
    \FORALL{solutions (sites) $s \in S$}{
      \STATE Evaluate the fitness (objective) of $s$
      \STATE Separate a Benders cut $c$
      \STATE Let $C := C \cup \{c\}$
    }\ENDFOR
    \REPEAT{
      \STATE Let $S = E\cup B\cup R$ where $|E|=ne$ and $|B|=nb$ with $R$ having lower fitness than the others.
      \FORALL{solutions (sites) $s \in E\cup B$}{
        \STATE Heuristically evaluate $nre$ (or $nrb$) solutions in the neighbourhood $ngh$ of $s$ if $s\in E$ (resp. $s\in B$)
        \IF{ better heuristic solution found}{
          \STATE Evaluate the fitness of the solution
          \STATE Separate a Benders cut $c$
          \STATE Let $C := C \cup \{c\}$
          \STATE Append the solution to $E$ (or $B$)
        }\ENDIF
      }\ENDFOR
      \STATE Let $S:=E\cup B \cup \{b\}$ where $b$ is the incumbent solution of the Benders decomposition
    }\UNTIL{stopping condition met}
  \end{algorithmic}
\end{algorithm}

\subsubsection{Initialization}
\label{subsubsec:initialization}

The algorithm is initialized with a population of $nre + nrb$ workers, which are
uniformly randomly distributed over the solution space. The fitness of each
solution is evaluated by solving the LP subproblem. Each LP subproblem produces
a Benders cut which is stored in a pool of cuts shared by the Benders
decomposition. The fitness scores are ranked and the $nb$ best ``flower
patches'' are selected for neighbourhood search. The algorithm enters the main loop.

Simultaneously, the algorithm commences solving the Benders decomposition
using the ``single tree'' master problem approach: Lazy constraint callbacks are
used to solve the LP subproblem and separate the cuts. This means that the
master problem need only be solved to optimality once as opposed to once per
iteration. Any generated cuts are add to the shared pool of cuts.

\subsubsection{The main loop}
\label{subsubsec:main_loop}

The main loop consists of two main phases: neighbourhood search and cut sharing.
The neighbourhood search is carried out by each process independently, while the
cut sharing represents a communication or synchronisation step between the
processes. 

\subsubsection{Neighbourhood search}
\label{subsubsec:neighbourhood_search}

Each iteration, the workers that discovered the $ne$ elite solutions each
recruit $nre$ workers for neighbourhood search. Likewise, the workers who
discovered the remaining $nb - ne$ best solutions each recruit $nrb$ workers
for neighbourhood search. The $ns - nb$ workers who failed to find a best solution
rejoin the pool of workers.

Neighbourhood search at a given site is performed by each worker producing a
pool of candidate solutions using a Hamming distance function which randomly
selects at most $ngh$ binary variables to alter. When solving arbitrary MIPs
with binary variables $x_i$ for $i\in N$ this is equivalent to imposing the
constraint $$ \sum_{i\in Z} x_i + \sum_{i\in N\setminus Z} (1-x_i) \le ngh.$$
That is, of all the variables $x_i,\ i\in Z$ which have $x_i=0$ in the current
solution and all the remaining variables that have $x_i=1$, only $ngh$ may
change their value. In our application a more specialised neighbourhood move can
be defined based on the structure of the problem. We have binary variables that
represent a unit increment in transmission capacity between two locations (a
right of way).  Each right of way (analogous to a set of edges $ij$ on a
multigraph with nodes $i$ and $j$) has $p$ binary variables denoting the
installation of a equivalent line. This means that individually installing the
$1^{st}$ line is equivalent to installing the $2^{nd}\dots p^{th}$ line. Clearly
it is undesirable for the Hamming distance function to randomly replace the
installation of one line on a right of way with another. For this reason the
function operates on groups of binary variables representing a single right of
way. If two or more changes are made to a right of way they are directionally
consistent i.e. if the first change added a line, subsequent changes will also
add a line until the maximum of $p$ lines are installed.

The fitness of each candidate solution in the pool is estimated using the
matheuristic given in Algorithm~\ref{alg:BBHA_heuristic}. Here the cost of the
master problem is calculated from the candidate solution. If the cost exceeds
the incumbent fitness value the evaluation stops. Otherwise, the shared set of
Benders cuts are used to estimate the cost of the subproblem. The piecewise
linear approximation of the original problem objective is of the form $c^T
y+\max_{i\in Cuts} \{r_i - B_i y\}$. Here $y$ are the master problem (binary)
decision variables with cost vector $c$. Each Benders cut takes the form $\hat
v\ge r_i - B_iy$ where $B_i$ is a row vector and $r_i$ a constant.  Collectively
these cuts enforce that $\hat v\ge  \max_{i\in Cuts} \{r_i - B_i y\}$ with
minimisation of $\hat v$ ensuring that at least one of the inequalities holds
with equality. Hence if $v=r-By$, the estimated fitness of the candidate
solution comprises the cost of the master problem plus the maximum value of the
vector $v$. 

\begin{algorithm}[htb!]
  \caption{Heuristic fitness evaluation}
  \begin{algorithmic}[1]{} 
  \label{alg:BBHA_heuristic}
    \REQUIRE{ $c:$ the vector of costs}
    \REQUIRE{ $\hat{y}:$ the candidate solution vector}
    \REQUIRE{ $B:$ the matrix of coefficients of known Benders cuts}
    \REQUIRE{ $r:$ the rhs vector of known Benders cuts}
    \STATE Let $C := c^T\hat{y}$
    \IF{$C\ \le\ $  fitness of the current incumbent}{
      \STATE Let $v := r - B\hat{y}$
      \STATE Let $C := C + max(v)$
    }\ENDIF
    \RETURN $C$
  \end{algorithmic}
\end{algorithm}


Each worker then solves the LP subproblem for the most promising heuristically
determined solution in their solution pool, and the generated Benders cuts are
added to the shared pool of known cuts.

The fitness scores of the solutions found by the workers are combined with the
incumbent solution of the Benders decomposition and are ranked from best to
worst. The $nb$ best solutions are selected for neighbourhood search during the
next iteration.

\subsubsection{Cut sharing}
\label{subsubsec:cut_sharing}

At the conclusion of the neighbourhood search phase any Benders cuts produced by
the Benders decomposition are added to the shared pool of cuts. Any cuts produced by the
workers are likewise made available to the Benders decomposition, and may be added
to the pool of cuts managed by CPLEX during a subsequent execution of the lazy
constraint callback. The effect of this cut sharing is that each of the bees has
a more accurate approximation of the objective function. With this approximation
the scout bee (branch and bound tree) avoids searching any solution that is not
at least potentially better than the best found so far. The worker bees use the
approximation to quickly evaluate solutions in the site to identify the most
promising solution for which a full LP is solved.

%
%
%

\subsubsection{Termination}
\label{subsubsec:termination}

The algorithm may terminate in several ways: After $n^{max}$ iterations,
$t^{max}$ seconds, or if the Benders decomposition identifies and proves the
optimal solution. Note that even if the time or iteration limit does not allow a
provably optimal solution to be found, we are still able to extract a lower
bound from the branch and bound tree that the scout bee was searching. Thus even
in the case where the method is only a heuristic, we have an estimate of the
maximum gap to the globally optimal solution.

\section{Mathematical model of TESP}
\label{sec:model}

We consider a electricity network consisting of nodes and arcs (referred to as
``rights of way''). The ability to generate power and demand for power occurs at
the nodes over a number of time periods. The objective of the complete TESP
model is the is to minimize the investment
cost of expanding the transmission network while simultaneously minimizing a
penalty for load curtailment at nodes with net demand. A discrete number of new
or reinforcing circuits may be installed on each right of way, and the
location and size of any ESS at the nodes are determined.

Cyclic discrete time is used to model the period of operation, and therefore the
state of any installed ESS in the last time interval must be identical to the
state in the initial time interval. This might model the typical power use over
a day with the end of one day matching the start of the next. Generation is
re-dispatchable and demand may vary between time intervals. Despite the
introduction of time to the model, the planning is static, and only a single
final expansion plan is produced. More complex models, for example with multiple
scenarios for demand and renewable power generation capacity, are possible but
not considered in this paper.

The model determines the network expansion plan, and operational characteristics
such as the amount of energy stored in the ESS, the network flows, and the phase
angles at each bus for each time interval. As with other variants of the
disjunctive TEP, power flows are modeled using a DC approximation
\cite[p.36]{knight197228}.

The mathematical model presented here, as well as alternative modelling
approaches in the literature, is discussed in detail in
\citep{macrae_benders_2016}. As such, only an abridged discussion of the
decomposed model follows. The key point to note here is that the main integer
(binary) variables to be determined relate to the rights of way to be installed,
while a very large number of continuous variables and associated constraints
have to be considered to determine the optimal generation, storage and power
flows for any choice of network expansion.

The following notation will be used throughout this paper to define the TESP:

\noindent \begin{flushleft}
\textit{Sets}\\[1ex]
\begin{tabular}{ll}
  $\Gamma$ & the set of indices for buses; \tabularnewline
  $\Omega_{0}$ & the set of rights of way for existing circuits; \tabularnewline
  $\Omega_{c}$ & the set of rights of way for candidate circuits; \tabularnewline
  $\Psi$ & the set of uniform time intervals $\left\{ 1,2,\ldots,T\right\} $;  \tabularnewline
\end{tabular}
\par\end{flushleft}

\noindent \begin{flushleft}
\textit{Parameters}\\[1ex]
\begin{tabular}{ll}
  $\alpha_{tk}$ & cost of curtailment at time $t$ at bus $k$; \tabularnewline
  $b_{k}$ & cost of installing storage at bus $k$; \tabularnewline
  $c_{ij}$ & cost of installing a circuit on right of way $ij$; \tabularnewline
  $d_{tk}$ & demand at time $t$ at bus $k$; \tabularnewline
  $\bar{f}_{ij}$ & maximum possible power flow on right of way $ij$; \tabularnewline
  $\bar{g}_{k}$ & maximum possible generation at bus $k$; \tabularnewline
  $\gamma_{ij}$ & susceptance of circuits installed on right of way $ij$; \tabularnewline
  $M_{ij}$ & the disjunctive parameter for right of way $ij$ \tabularnewline
  $n^{0}_{ij}$ & number of existing circuits on right of way $ij$; \tabularnewline
  $\bar{n}_{ij}$ & maximum number of installable circuits on right of way $ij$; \tabularnewline
  $\bar{x}_{k}$ & maximum installable storage capacity at bus $k$; \tabularnewline
  $\hat{y}^{p}_{ij}$ & binary parameter denoting installation of the $p$\textsuperscript{th} candidate circuit on right of way $ij$; \tabularnewline
\end{tabular}
\end{flushleft}

\noindent\begin{flushleft}
  \textit{Decision variables}\\[1ex]
\begin{tabular}{ll}
  $\beta_{tk}$ & power flow to storage at bus $k$ at time~$t$; \tabularnewline
  $g_{tk}$ & generation at time $t$ at bus $k$; \tabularnewline
  $f^{0}_{tij}$ & power flow for existing circuits at time $t$ on right of way $ij$; \tabularnewline
  $f^{p}_{tij}$ & power flow for the $p$\textsuperscript{th} candidate circuit at time $t$ on right of way $ij$; \tabularnewline
  $l_{tk}$ & level of storage at bus $k$ at time $t$;  \tabularnewline
  $r_{tk}$ & demand curtailment at time $t$ at bus $k$ ; \tabularnewline
  $\theta_{tk}$ & phase angle at time $t$ at bus $k$; \tabularnewline
  $x_{k}$ & storage capacity installed at bus $k$; \tabularnewline
  $y^{p}_{ij}$ & binary variable denoting installation of the $p$\textsuperscript{th} candidate circuit on right of way $ij$; \tabularnewline
  $v$ & estimate of the contribution of the subproblem to the objective function of the master problem \tabularnewline
\end{tabular}
\end{flushleft}

\subsection{The master problem}
\label{subsec:master}

\abovedisplayskip=-1ex
\belowdisplayskip=1ex
The objective of the master problem is to minimize the function

\begin{equation}
z = \underset{(i,j)}{\sum}c_{ij}y_{ij}^{p} + v
\label{eq:reform_master_obj}
\end{equation}
where $c_{ij}$ is cost of installing a line on right of way $ij$ and $y_{ij}^{p}$
is a binary variable denoting the installation of the $p$\textsuperscript{th}
candidate line on $ij$. The estimated contribution of the subproblem to the
objective function is given by $v$.

The following constraints are necessary to the master problem:
\\[-0.5ex]

\noindent \textit{Symmetry breaking constraints}

\begin{equation}
  y^{p}_{ij} \geq y^{p+1}_{ij} \quad\forall\;
  \left(i,j\right)\in\Omega_{c},\;\forall\;p\in
  \left\{1\ldots\bar{n}_{ij}-1\right\} \label{eq:reform_symmetry}
\end{equation}

The lexicographical constraint (\ref{eq:reform_symmetry}) eliminates the
symmetry of the binary decision variables by mandating the order of installation
of parallel circuits be arbitrary.
\\[-0.5ex]

\noindent\textit{Other constraints}
\begin{eqnarray}
  v &\geq &0  \label{eq:reform_master_bounds}\\
  y_{ij}^{p}&\in&\mathrm{\{0,1\}}
  \label{eq:reform_master_binary}
\end{eqnarray}

\subsection{The subproblem}
\label{subsec:subproblem}

Given a set of new circuit installations determined by the master
problem, the subproblem determines the cost of any installed ESS, and a
penalty for load curtailment.

\abovedisplayskip=-1ex
\belowdisplayskip=1ex
The objective of the subproblem is to minimize the function

\begin{equation}
v =
\sum_{k\in\Gamma}b_{k}x_{k} +
\sum_{t \in \Psi} \sum_{k\in\Gamma} \alpha_{tk} r_{tk}
\end{equation}
where $b_{k}$ is the fixed cost of installing $x_{k}$ MW of storage at bus $k$,
and $\alpha{tk}$ the cost of curtailing $r_{tk}$ in each time interval $t$. It is
assumed that the variable operating cost of ESS is low relative to fixed costs,
and these are therefore omitted from the objective function.

The following technical constraints govern the operation of the network:
\\[-0.5ex]

\noindent\textit{Nodal balance and power flow}

\begin{equation}
  \zeta + g_{tk} + r_{tk} - \beta_{tk}=d_{tk} \quad \forall\; t\in\Psi,\;\forall\; k\in\Gamma
   \label{eq:KCL}
\end{equation}
where
\begin{equation}
 \zeta = \underset{\left(i,k\right) \in \Omega_{0}}{\sum}f^{0}_{tik}
  - \underset{\left(k,j\right) \in \Omega_{0}}{\sum}f^{0}_{tkj}
  + \overset{\bar{n}_{ij}}{\underset{p=1}{\sum}}
    \underset{\left(i,k\right) \in \Omega_{c}}{\sum}f^{p}_{tik}
    - \overset{\bar{n}_{ij}}{\underset{p=1}{\sum}}
    \underset{\left(k,j\right) \in \Omega_{c}}{\sum}f^{p}_{tkj}
  \end{equation}
Nodal balance i.e. Kirchhoff's current law is ensured for each time interval by
constraint (\ref{eq:KCL}).

Power flows are modeled using a DC approximation requiring that the phase angle
at each bus be determined for each time interval:

\begin{equation}
  \begin{split}
    f_{tij}^{0}-\gamma_{ij}n_{ij}^{0}\left(\theta_{ti}-\theta_{tj}\right)=0
    \quad
    \forall\; t\in\Psi,\forall\;\left(i,j\right)\in\Omega_{0} \label{eq:KVL_existing}
  \end{split}
\end{equation}
\begin{equation}
  \begin{split}
    \lvert f_{tij}^{p}-\gamma_{ij} \left(\theta_{ti} -
    \theta_{tj}\right)\rvert  \leq M_{ij}(1-\hat{y}_{ij}^{p}) \quad
  \forall\; t\in\Psi,\;\forall\;\left(i,j\right)\in\Omega_{c},\;\forall\;p\in \left\{1\ldots\bar{n}_{ij}\right\} \label{eq:KVL_candidate}
  \end{split}
\end{equation}

Kirchhoff's voltage law is implemented for existing circuits by
(\ref{eq:KVL_existing}), and for candidate circuits by (\ref{eq:KVL_candidate}).
Note that absolute values are given to simplify the notation, however in
practice these are
readily expanded into pairs of ranged linear constraints.
The disjunctive parameter $M_{ij}$ must be large enough that it does not limit
the difference in phase angles of buses $i$ and $j$. Minimal values of $M_{ij}$
may be calculated by following the procedure given in \citep{binato_new_2001}.
\begin{eqnarray}
  \left\lvert f_{tij}^{0}\right\rvert&\leq& n_{ij}^{0}\bar{f}_{ij}\quad\forall\; t\in\Psi,\;\forall\;\left(i,j\right)\in\Omega_{0}
\label{eq:powerflow_existing}\\
    \left\lvert f_{tij}^{p}\right\rvert &\leq &\hat{y}_{ij}^{p}\bar{f}_{ij}\quad\forall\; t\in\Psi,\;\forall\;\left(i,j\right)\in\Omega_{c},\;\forall\;p\in \left\{1\ldots\bar{n}_{ij}\right\}
\label{eq:powerflow_candidate}
\end{eqnarray}

Constraint (\ref{eq:powerflow_existing}) and constraint
(\ref{eq:powerflow_candidate}) enforce nominal thermal limits on existing and
candidate circuits respectively.
\\[-0.5ex]

\noindent \textit{Storage level and charge/discharge limits}

\begin{equation}
l_{1k}=l_{Tk}+\beta_{1k}\quad\forall\; k\in\Gamma
\label{eq:level_start}
\end{equation}
\begin{equation}
l_{tk}=l_{t-1,k}+\beta_{tk}\quad\forall\; t\in\Psi,\;\forall\; k\in\Gamma
\label{eq:level_mid}
\end{equation}

The set of time intervals $\Psi$ is assumed to be cyclic to allow the operation of
the storage throughout the desired time period, for example, a typical day. As
such, the storage level at the end of the day is required to match the initial
storage state. The ``wrap around'' constraint (\ref{eq:level_start}) implements
this requirement. For all other time intervals the storage level is determined by
(\ref{eq:level_mid}).

\begin{equation}
0\leq l_{tk}\leq x_{k}\quad\forall\; t\in\Psi,\;\forall\; k\in\Gamma \label{eq:level_bounds}
\end{equation}
\begin{equation}
  0 \leq x_{k} \leq \bar{x}_{k}\quad\forall\; k\in\Gamma
  \label{eq:storage_bounds}
\end{equation}
Constraint (\ref{eq:storage_bounds}) establishes bounds on the installable
storage capacity at bus $k$, while constraint (\ref{eq:level_bounds}) ensures the
stored energy does not exceed the installed capacity.
\\[-0.5ex]

\noindent\textit{Generation bounds}

\begin{equation}
0\leq g_{tk}\leq\bar{g}_{k}\quad\forall\; t\in\Psi,\;\forall\; k\in\Gamma
\label{eq:generation_bounds}
\end{equation}
Constraint (\ref{eq:generation_bounds}) imposes bounds on generator re-dispatch.
\\[-0.5ex]

\noindent \textit{Curtailment bounds}

\begin{equation}
  0\leq r_{tk} \leq d_{tk}\quad\forall\; t\in\Psi,\;\forall\; k\in\Gamma \label{eq:curtailment_bounds}
\end{equation}
Load curtailment at any bus $k$ cannot exceed demand during the same time
interval $t$..
\\[-0.5ex]

\noindent\textit{Other constraints}

\begin{equation}
  f^{0}_{tij}, f^{p}_{tij}, \beta_{tk}, \theta_{tk}\:\mathrm{unbounded}
\end{equation}

\subsection{Optimality cut}
\label{subsec:opt}

As noted above, load curtailment is permitted at any bus during any time
interval so long as it does not exceed demand at that bus during the same time
period. Therefore, the dual of the subproblem remains bounded for any feasible
solution to the master problem. Accordingly, we need only consider the following
optimality cut:

Let the dual variables $\pi_{d_{tk}}$ correspond to constraint
(\ref{eq:KCL}), $\pi_{\gamma_{tij}}$ to constraint (\ref{eq:KVL_existing}),
$\pi_{\gamma^{+p}_{tij}}$ and $\pi_{\gamma^{-p}_{tij}}$ to
constraint (\ref{eq:KVL_candidate}), $\pi_{f^{+0}_{tij}}$ and
$\pi_{f^{-0}_{tij}}$ to (\ref{eq:powerflow_existing}), and $\pi_{f^{+p}_{tij}}$
and  $\pi_{f^{-p}_{tij}}$ to (\ref{eq:powerflow_candidate}). The dual
variables $\pi_{s_{tk}}$ correspond to constraints(\ref{eq:level_start})
and (\ref{eq:level_mid}), and $\pi_{\bar{l}_{k}}$ to (\ref{eq:level_bounds}).
Lastly, let the dual variables $\pi_{g_{tk}}$, $\pi_{r_{tk}}$, and $\pi_{x{k}}$
correspond to the bounds (\ref{eq:storage_bounds} - \ref{eq:curtailment_bounds})
respectively.

The optimality cut is therefore

\begin{equation}
  \begin{split}
  v - \underset{t \in \Psi} {\sum}
    \underset{\left(i,j\right) \in \Omega_{c}} {\sum}
    \left[ \pi_{f^{+p}_{tij}}y_{ij}^{p}\bar{f_{ij}}
    + \pi_{f^{-p}_{tij}}y_{ij}^{p}\bar{f_{ij}} \right]
  - \underset{t \in \Psi} {\sum}
    \underset{\left(i,j\right) \in \Omega_{c}} {\sum}
    \left( \pi_{\gamma^{+p}_{tij}} + \pi_{\gamma^{-p}_{tij}}\right)
    \left(M_{ij}(1-y_{ij}^{p})\right)
    \geq \\
    \underset{t \in \Psi} {\sum}
    \underset{k \in \Gamma} {\sum} d_{tk}\pi_{d_{tk}}
  + \underset{t \in \Psi} {\sum}
    \underset{\left(i,j\right) \in \Omega_{0}} {\sum}
    \left[ \pi_{f^{+0}_{tij}}n_{ij}^{0}\bar{f_{ij}}
    + \pi_{f^{-0}_{tij}}n_{ij}^{0}\bar{f_{ij}} \right]
  + \underset{t \in \Psi} {\sum}
    \underset{k \in \Gamma} {\sum} \bar{g_{k}}\pi_{g_{tk}}
  + \underset{t \in \Psi} {\sum}
    \underset{k \in \Gamma} {\sum} d_{tk}\pi_{r_{tk}}
  + \underset{k \in \Gamma} {\sum} \bar{x_{k}}\pi_{x_{k}}
    \label{eq:sd_obj}
  \end{split}
\end{equation}

\subsection{Limitations}
\label{subsec:limitations}

The relative simplicity of this TESP formulation comes at the cost of addressing
certain features of a real world electrical transmission network.
The most obvious limitation is that power flows are modelled using a DC
approximation to the AC power flow of most transmission networks.

It is also assumed that the variable operating cost of ESS is negligible, at least
relative to the fixed cost of installing and operating the ESS over its
lifetime, and that fixed costs increase linearly with capacity. Furthermore,
power flow to and from ESS is limited only be the total capacity and current
level of the storage. The model allows that the storage completely charge or
discharge within a single time interval. Furthermore, the model assumes 100\%
efficiency for storage and losses are not considered.

The model also assumes generator re-dispatch does not incur cost, and that
generators are not subject to technical constraints such as generation
ramp rates. 

While it is possible to address these and other limitations of the model with
additional variables and constraints, these come at the cost of significant
complexity in both notation and implementation. Here we have sought to balance
to the realism of the modelling with the intent to use the model simply to
demonstrate the use of the algorithmic approach.

\section{Numerical results}
\label{sec:num_res}

In each of the numerical experiments described in this section the model is
implemented in Python 3.4.3 and, where appropriate, makes use of the Python
library for IBM ILOG CPLEX 12.6.3. Parallelization is achieved using multiple
processes, not threading. The Benders decomposition is implemented with a
``single tree'' master using lazy constraint callbacks. Preprocessing is
disabled by default, and while the LP solver may take advantage of
multi-threading, the branch and cut is single threaded.

\subsection{Parameter tuning}
\label{subsec:param_tuning}

There are a number of parameters to the BBHA algorithm which may be tuned to
find a set of default values that empirically demonstrate good performance. The
tuneable parameters are given in Table~\ref{tab:params}.

The IEEE-25 bus test system is used to benchmark combinations of parameters
presumed likely to perform well. A schematic and tabulated data are available in
\citep{ekwue_transmission_1984}. The system has 25 buses and 36 rights of way
with a total demand of 2750 MW. Without storage, and permitting a maximum of 4
new or reinforcing circuits on each right of way, the cost of the optimal
expansion plan is US\$107.7 million.

While it would be preferable to incorporate real world storage costs into the
model, the cost per MW of long term energy storage technology is currently
high enough to prevent the installation of any storage in the test systems
discussed in this paper. Therefore, an arbitrary cost coefficient of
US\$2000MW/h is used for each network to ensure storage is installed.

\begin{figure}[!ht]
  \centering
  \includegraphics[width=.75\columnwidth]{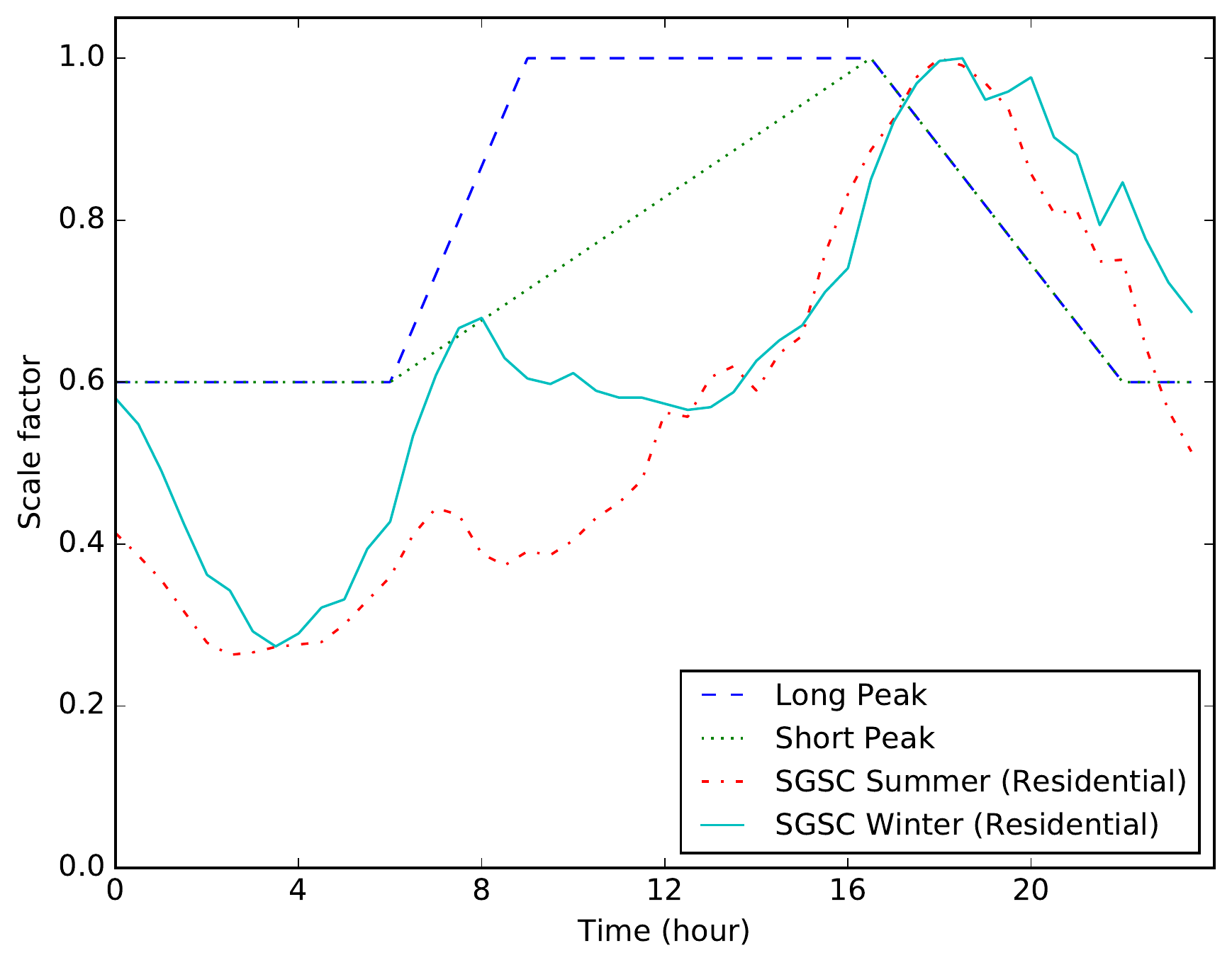}
  \caption{Load profiles used for each case study. (24 hours at a 30 minute
  interval).}
  \label{fig:load_profiles}
\end{figure}

Under the long peak scenario shown in Figure~\ref{fig:load_profiles} the cost
of the optimal expansion plan is US\$43.8 million. This result is the
benchmark objective for the parameter tuning.

\begin{figure}[!ht]
  \centering
  \includegraphics[width=.75\columnwidth]{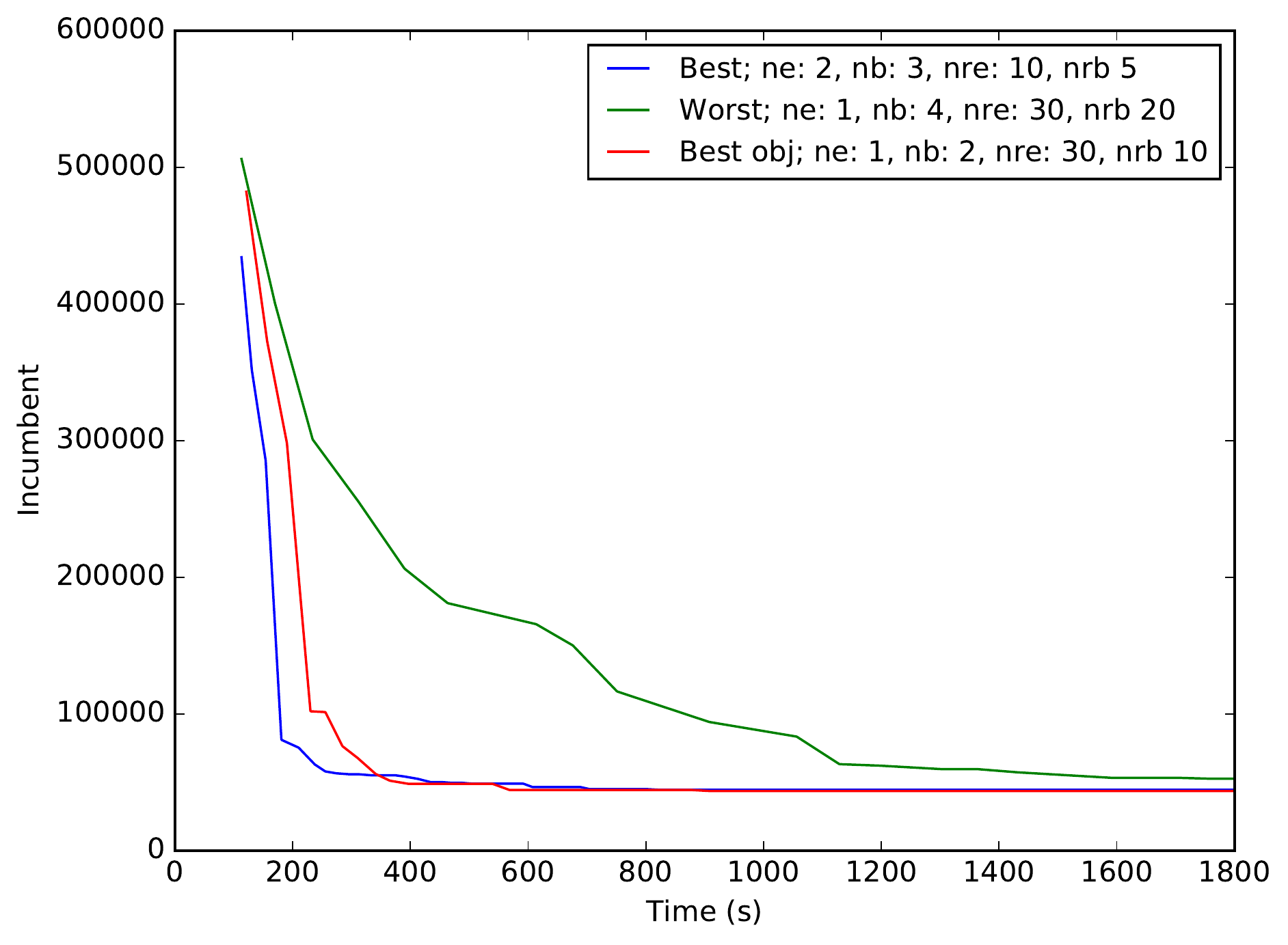}
  \caption{Best and worst parameter sets for IEEE 25-bus network and long peak
    scenario. The ``Best'' produces the best solution quality on average over the run,
  while ``Best obj'' produces the optimal solution the fastest for this instance.}
  \label{fig:tuning}
\end{figure}

In this tuning exercise, 34 sets of parameters are compared over the first 1800
seconds (30 minutes) of the optimization. The results are given in
Table~\ref{tab:tuning}. We use a composite trapezoidal rule to integrate along
the time axis and then rescale against the worst (largest) integral (scaled
trapz). The parameter sets are ranked then ranked. The best result is that
with lowest value. Only 5 sets of parameters ([ne: 1, nb: 2, nre: 30, nb: 10],
[ne: 1, nb: 2, nre: 30, nb: 15], [ne: 2, nb: 3, nre: 20, nb: 10], [ne: 3, nb:
4, nre: 10, nb: 5], and [ne: 3, nb: 4, nre: 20, nb: 15]) find the optimal
solution within the 30 minute window. The timeseries of the incumbent value of
the best of these parameter sets is plotted with the best and worst parameter
sets in Figure~\ref{fig:tuning}.

In general, parameters sets with a relatively modest number of workers and
associated high number of iterations appear to do well. An exception is the
parameter set (ne: 1, nb: 4, nre: 20, nb: 15) which requires only 39
iterations to match the best sub-optimal objective function value. This is
explained by the proximity of the 3 elite search neighbourhoods, and subsequent
thorough exploration of the combined neighbourhood. Other similar parameter sets
that match this objective function value by the end of the optimization do not
converge as quickly, as evidenced by their larger scaled trapz scores.

The Benders scout ensures that the BBHA is guaranteed to find the exact optimal
solution to the problem given sufficient time to run to completion. Of course
this may take a significant amount of time. The objective of the BBHA is to
discover high quality solutions quickly, and as such we favour parameter sets
which rapidly converge to such solutions in the case studies that follow. We
explore the performance of three sets of parameters:


\begin{itemize}
  \item[] [ne: 2, nb: 3, nre: 10, nrb: 5]: the parameter set with the smallest scaled
    trapz measure.
  \item[] [ne: 1, nb: 2, nre: 30, nrb: 10]: the parameter set that converges to
    the optimal solution the fastest using the scaled trapz measure.
  \item[] [ne: 1, nb: 2, nre: 10, nrb: 5]: the parameter set with the largest
    number of iterations.
\end{itemize}

\begin{figure}[!ht]
  \centering
  \includegraphics[width=.75\columnwidth]{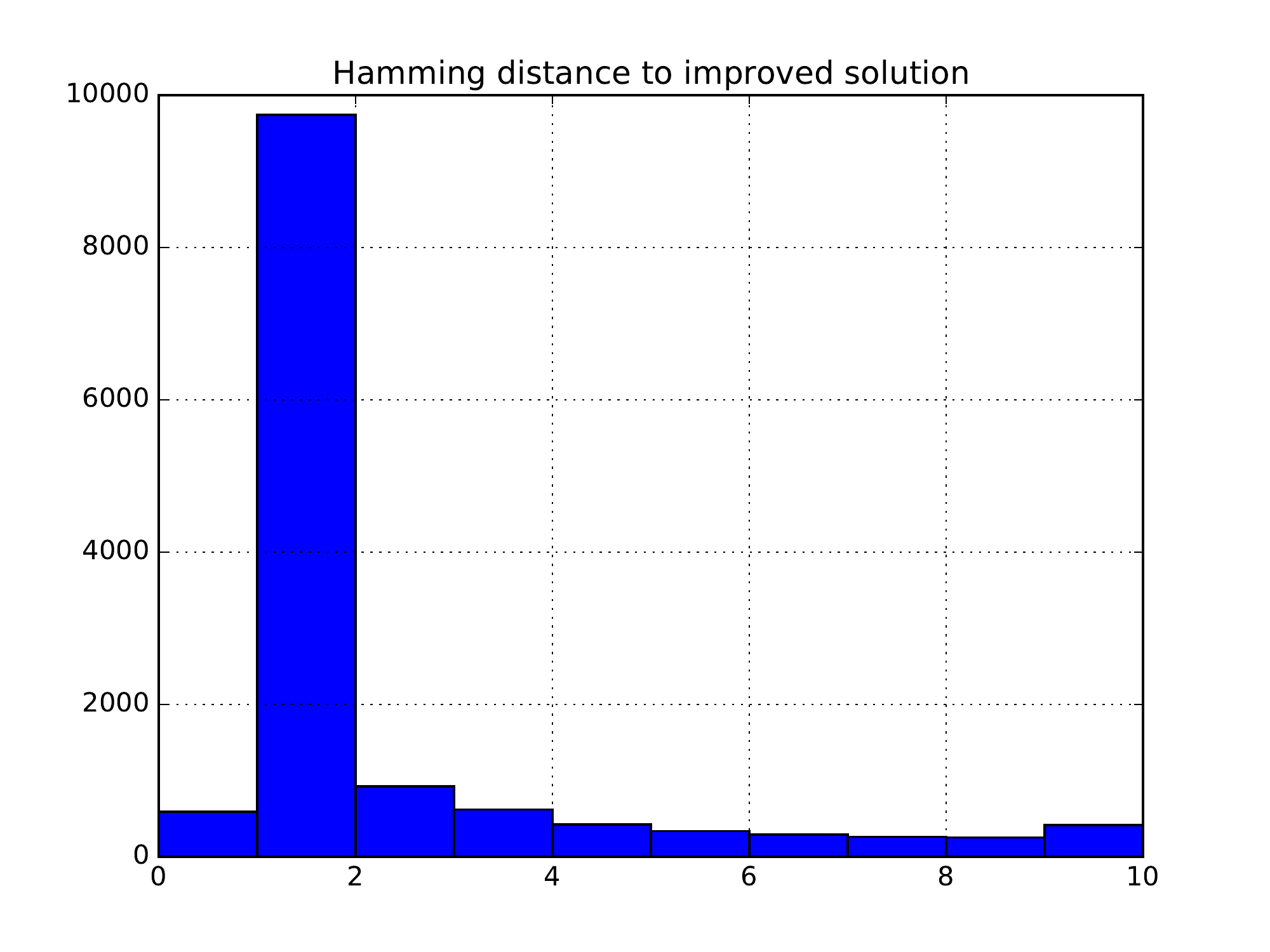}
  \caption{Histogram showing the hamming distance to best improved solution.}
  \label{fig:hamminghist}
\end{figure}

The final parameter to consider is the size of the neighbourhood for local
search $ngh$, described in Subsection \ref{subsubsec:neighbourhood_search}.
A histogram showing the distribution of the hamming distance over the range 1-10
required to produce the best improved solution of nearly 14000 workers is shown
in Figure~\ref{fig:hamminghist}. A value of 2 accounts for the largest number of
improved solutions. This is perhaps unsurprising as it reflects the somewhat
routine circuit swap in which one circuit is deselected and another selected.
Given that the long tail of larger hamming distances typically resulted in
improved solutions only at the beginning of the optimization the value of $ngh$
was reduced to 8 for the case studies.

\subsection{Case study: 46-bus network}
\label{subsec:cs_46-bus}

Representing the southern part of the Brazilian transmission network, the 46-bus
test system comprises 46 buses and 79 rights of way. Total demand in the network
is 6880MW. Tabulated data are provided in \citep{haffner_branch_2000}. The
investment cost of the optimal expansion plan without ESS is US\$154.42 million.

In this case study we allow the installation of a maximum of 5 new or
reinforcing circuits on each right of way. Storage may be installed at any bus
at an arbitrary cost of US\$2000MW/h.

As the amount of storage installed depends upon the demand scenario under which
it is operated, four demand scenarios are considered. The short peak and long
peak scenarios are described in \citep{macrae_benders_2016}, and the Smart Grid,
Smart City (SGSC) residential winter and summer scenarios are generic load
profiles taken from \citep{berry_national_2013}. Each scenario describes a 24
hour period with a 30 minute time step, and is shown in
Figure~\ref{fig:load_profiles}.

Each scenario is optimized $N=5$ times for both the BBHA and Bees algorithm, and
once using Benders decomposition which as a deterministic method exhibits little
variance. Each optimization is limited to 4 hours. Tabulated results are given
in Table~\ref{tab:results}.

\begin{figure}[!ht]
  \centering
  \includegraphics[width=.75\columnwidth]{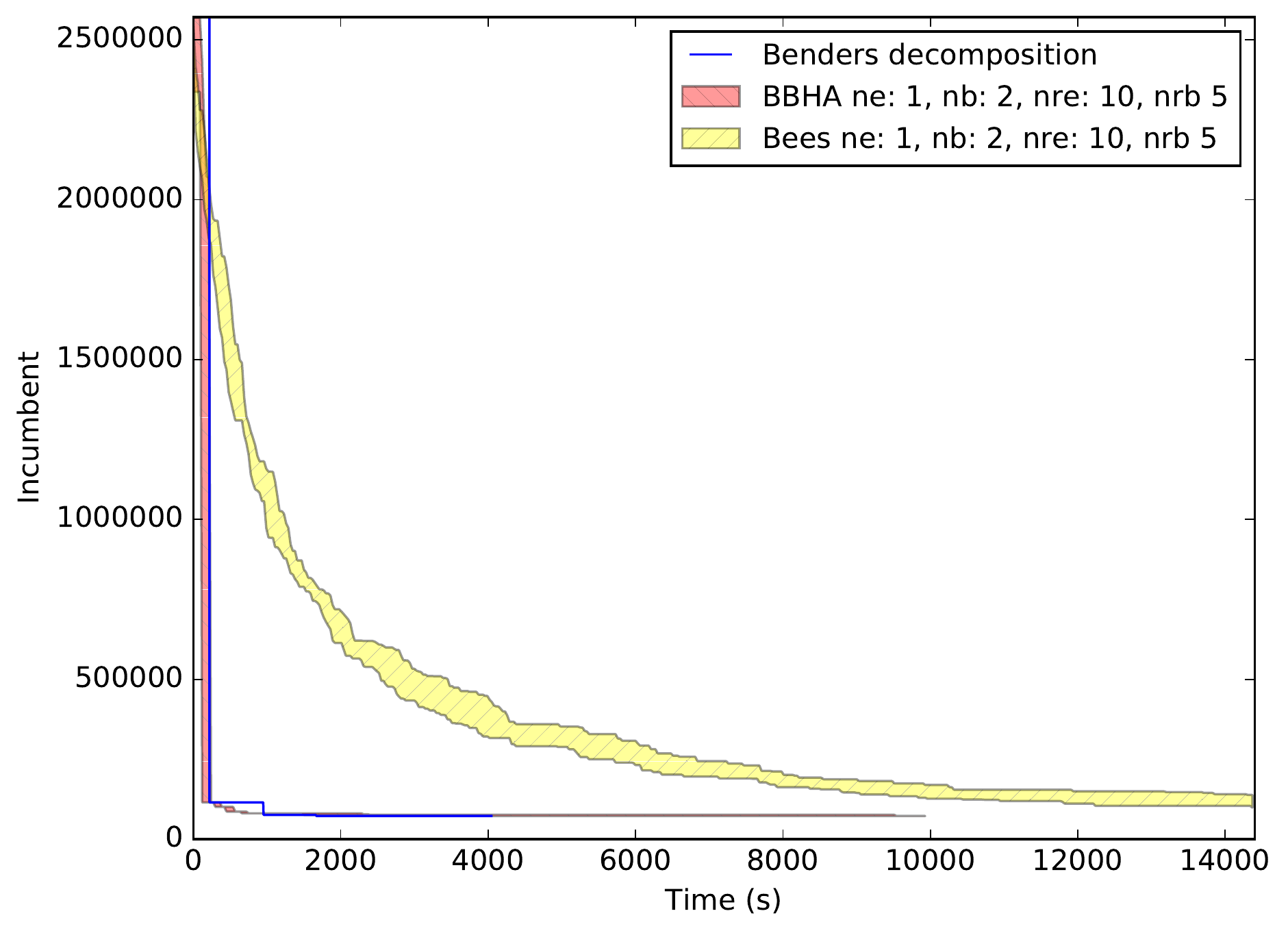}
  \caption{46-bus Short Peak scenario ne: 1, nb: 2, nre: 10, nrb: 5}
  \label{fig:46short_1_2_10_5}
\end{figure}

The parameter set [ne: 1, nb: 2, nre: 10, nrb: 5] typically matches or exceeds
the mean performance of the other parameter sets under investigation for the
BBHA, whereas the parameter set [ne: 1, nb: 2, nre: 30, nrb: 10]
exhibits better performance for the BA. The BBHA finds the optimal
solution for the short peak and SGSC summer and SGSC winter demand scenarios,
and the Benders scout is able to prove optimality. This is also true of the
Benders decomposition run. The range of incumbent solution values over time are
shown in Figures~\ref{fig:46short_1_2_10_5}, \ref{fig:46sum_1_2_10_5} \&
\ref{fig:46win_1_2_10_5} for the short peak, SGSC summer and SGSC winter
scenarios respectively.

\begin{figure}[!ht]
  \centering
  \includegraphics[width=.75\columnwidth]{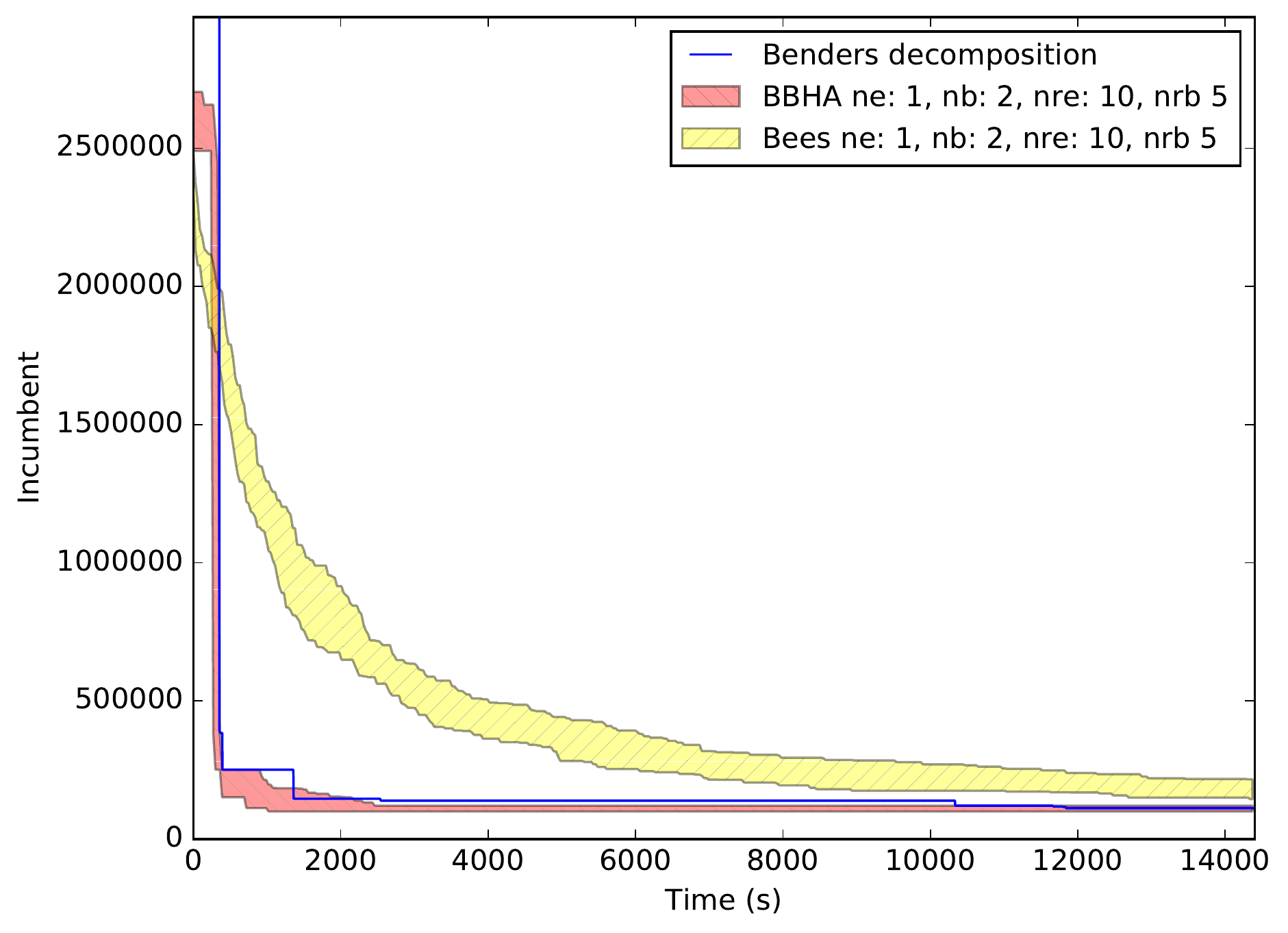}
  \caption{46-bus Long Peak scenario ne: 1, nb: 2, nre: 10, nrb: 5}
  \label{fig:46long_1_2_10_5}
\end{figure}

For the long peak scenario the best BBHA runs find the optimal solution, but
optimality is not proven. However, it is possible to use the pool of generated
cuts to prove a lower bound if necessary. A plot of the range of incumbent
solution values over the duration of the optimization is given in
Figure~\ref{fig:46long_1_2_10_5}.


\subsection{Case study: 93-bus network}
\label{subsec:cs_93-bus}

The Colombian 93-bus network is a medium complexity transmission network with 93
buses and 155 possible rights-of-way. The planning horizon includes 3 discrete
stages making this test system useful for testing multi-stage optimization
techniques \citep{escobar_multistage_2004}. In this case study will consider
only the total demand of 14559 MW in the final stage of the planning horizon.


A maximum of 4 new or reinforcing circuits is permitted to be installed on each
right of way. As with the previous case storage may be install at any bus at an
arbitrary cost of US\$2000MW/h. Network expansion plans are optimized for the
long peak, short peak, SGSC summer, and SGSC winter scenarios over a 4 hour
period. Tabulated results are included Table~\ref{tab:results}.

\begin{figure}[!ht]
  \centering
  \includegraphics[width=.75\columnwidth]{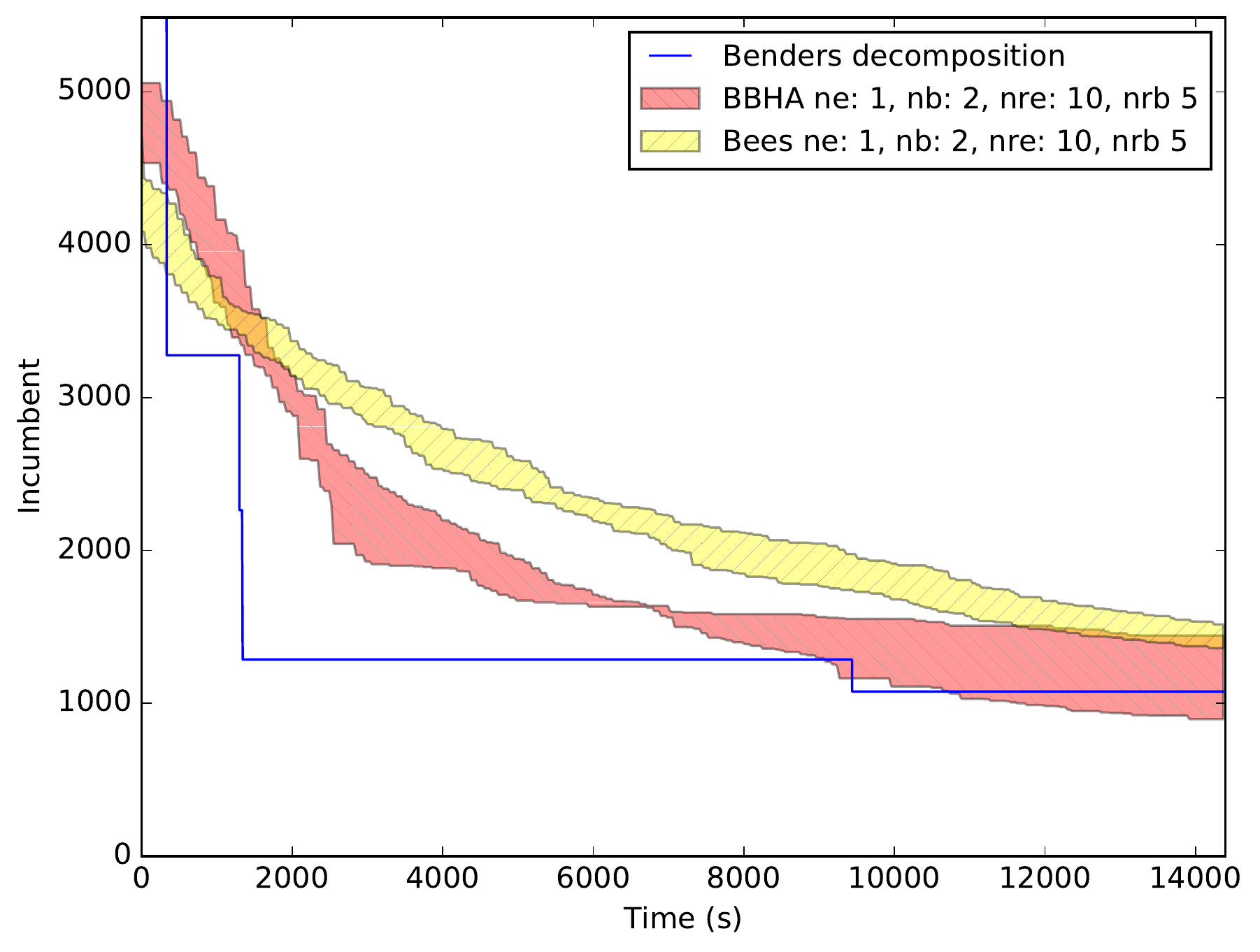}
  \caption{93-bus SGSC Winter scenario ne: 1, nb: 2, nre: 10, nrb: 5}
  \label{fig:93win_1_2_10_5}
\end{figure}

\begin{figure}[!ht]
  \centering
  \includegraphics[width=.75\columnwidth]{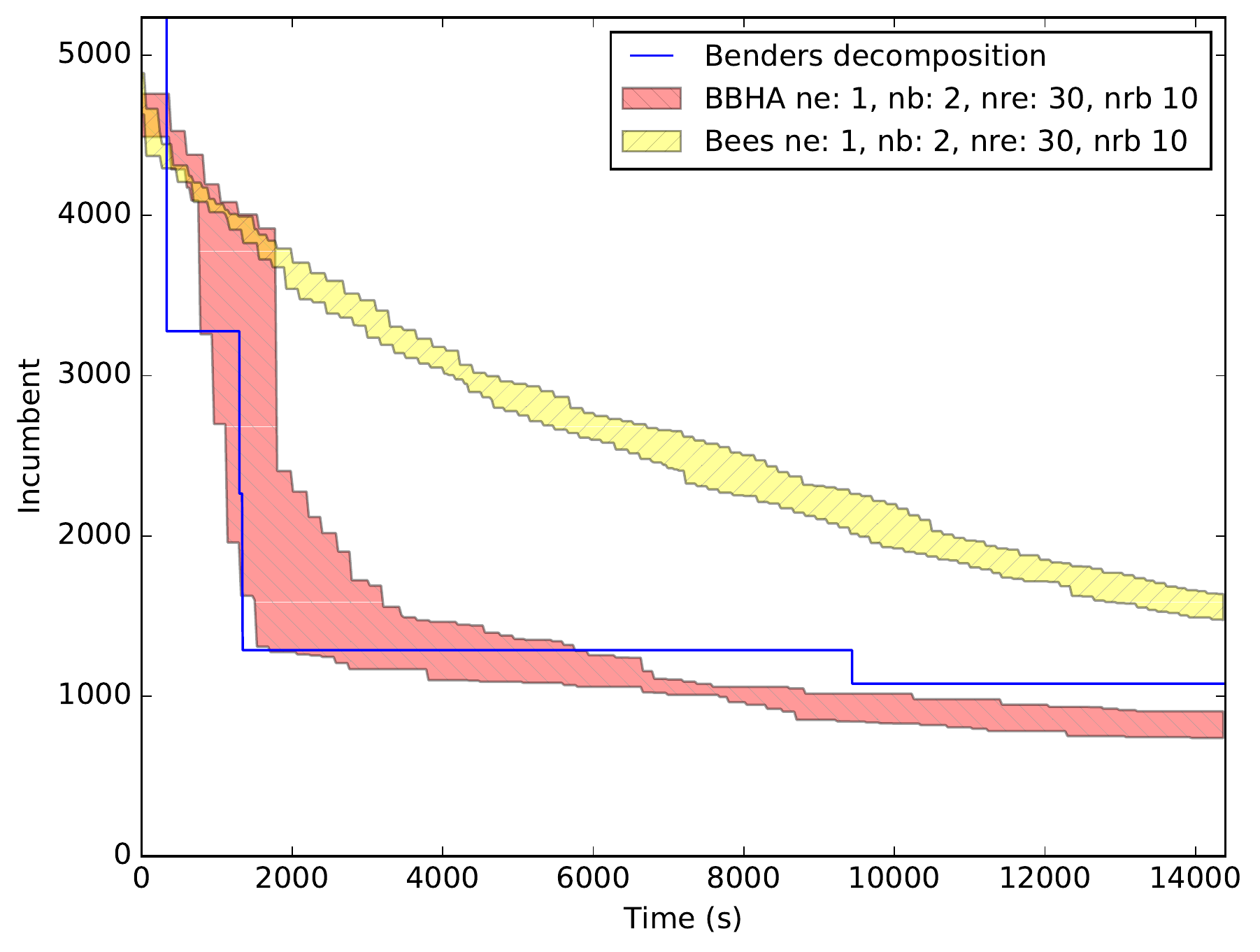}
  \caption{93-bus SGSC Winter scenario ne: 1, nb: 2, nre: 30, nrb: 10}
  \label{fig:93win_1_2_30_10}
\end{figure}

For this test system the parameter set [ne: 1, nb: 2, nre: 10, nrb: 5] exhibits
consistently good performance for both the BBHA and BA. The BA achieves a lower
mean for the long peak scenario. The Benders decomposition tends to lag behind
both approaches for all scenarios except the SGSC winter demand
profile, as shown in Figures~\ref{fig:93win_1_2_10_5} \& \ref{fig:93win_1_2_30_10}.

\subsection{Discussion}
\label{subsec:discussion}

The BBHA exhibits the essential characteristics of a hybrid optimization method.
Where the problem is readily solved by one of the component optimization methods
the BBHA performs comparably at minimum. Where each component optimization
method performs similarly on a given problem, the hybrid approach exceeds this
individual performance. In short, the whole is greater than the sum of its
parts.

\begin{figure}[!ht]
  \centering
  \includegraphics[width=.75\columnwidth]{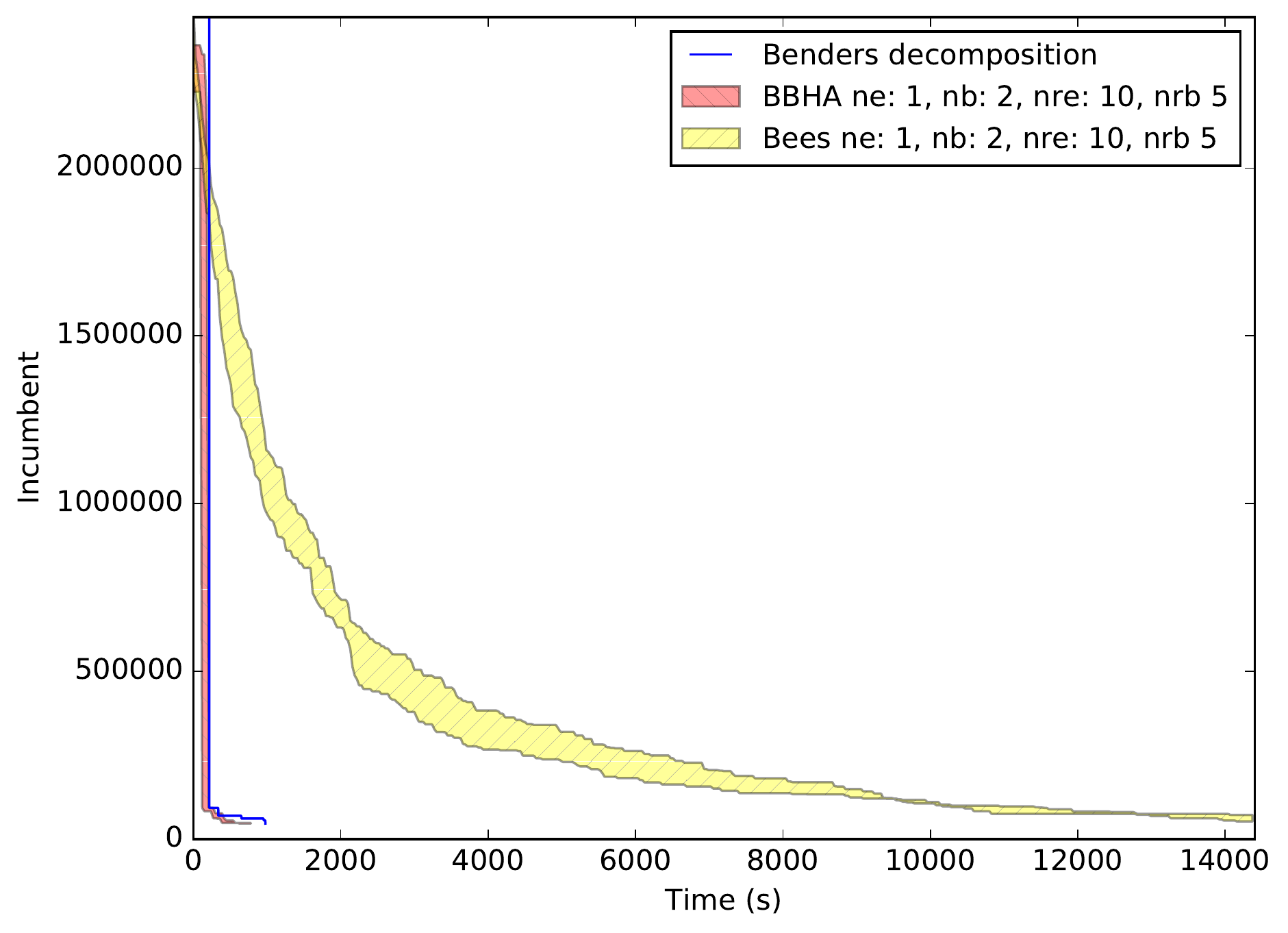}
  \caption{46-bus SGSC Summer scenario ne: 1, nb: 2, nre: 10, nrb: 5}
  \label{fig:46sum_1_2_10_5}
\end{figure}

\begin{figure}[!ht]
  \centering
  \includegraphics[width=.75\columnwidth]{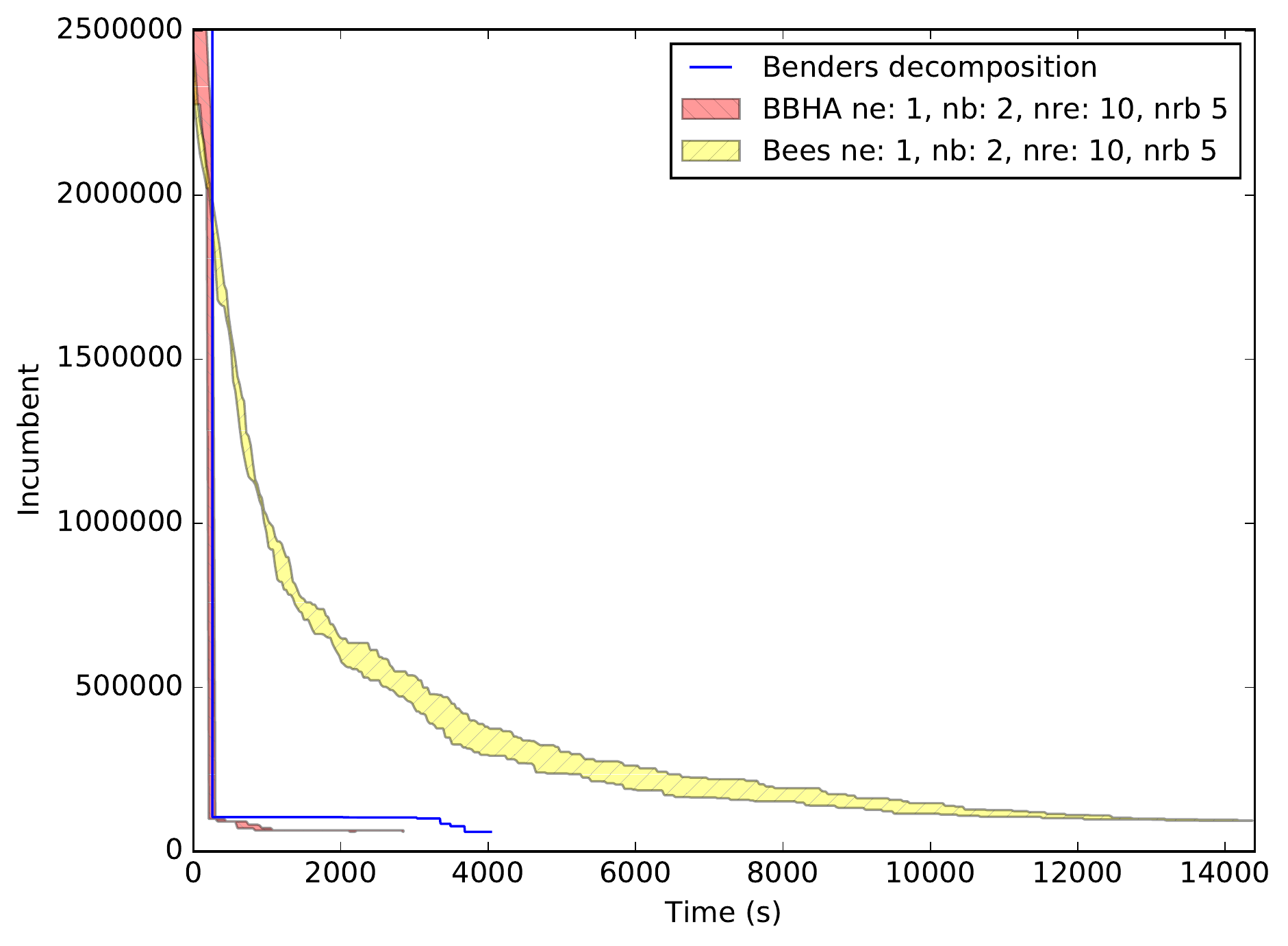}
  \caption{46-bus SGSC Winter scenario ne: 1, nb: 2, nre: 10, nrb: 5}
  \label{fig:46win_1_2_10_5}
\end{figure}

Figure~\ref{fig:46short_1_2_10_5} shows the 46-bus test system with the short
peak demand scenario, a problem known to amenable to Benders decomposition.
The BBHA performs comparably to the Benders decomposition, and in most runs
discovers the optimal solution earlier.  This can be observed by the incremental
improvements to the incumbent value over the first 1000 seconds. Both methods
are able to prove optimality within the time limit, however in this case the BBHA
takes longer. For the SGSC summer and SGSC winter scenarios shown in
Figures~\ref{fig:46sum_1_2_10_5} and \ref{fig:46win_1_2_10_5} respectively, the
BBHA not only discovers the optimal solution heuristically well in advance of
the Benders decomposition, but is also able to prove optimality prior.

In the case of the 46-bus test system and long peak scenario shown in
Figure~\ref{fig:46long_1_2_10_5}, the best BBHA run discovers the optimal
solution but is not able to prove it optimal within the time limit. The mean
solution is slightly better than the incumbent of the Benders decomposition,
however the worst BBHA solution is 109\% of the Benders decomposition incumbent.

\begin{figure}[!ht]
  \centering
  \includegraphics[width=.75\columnwidth]{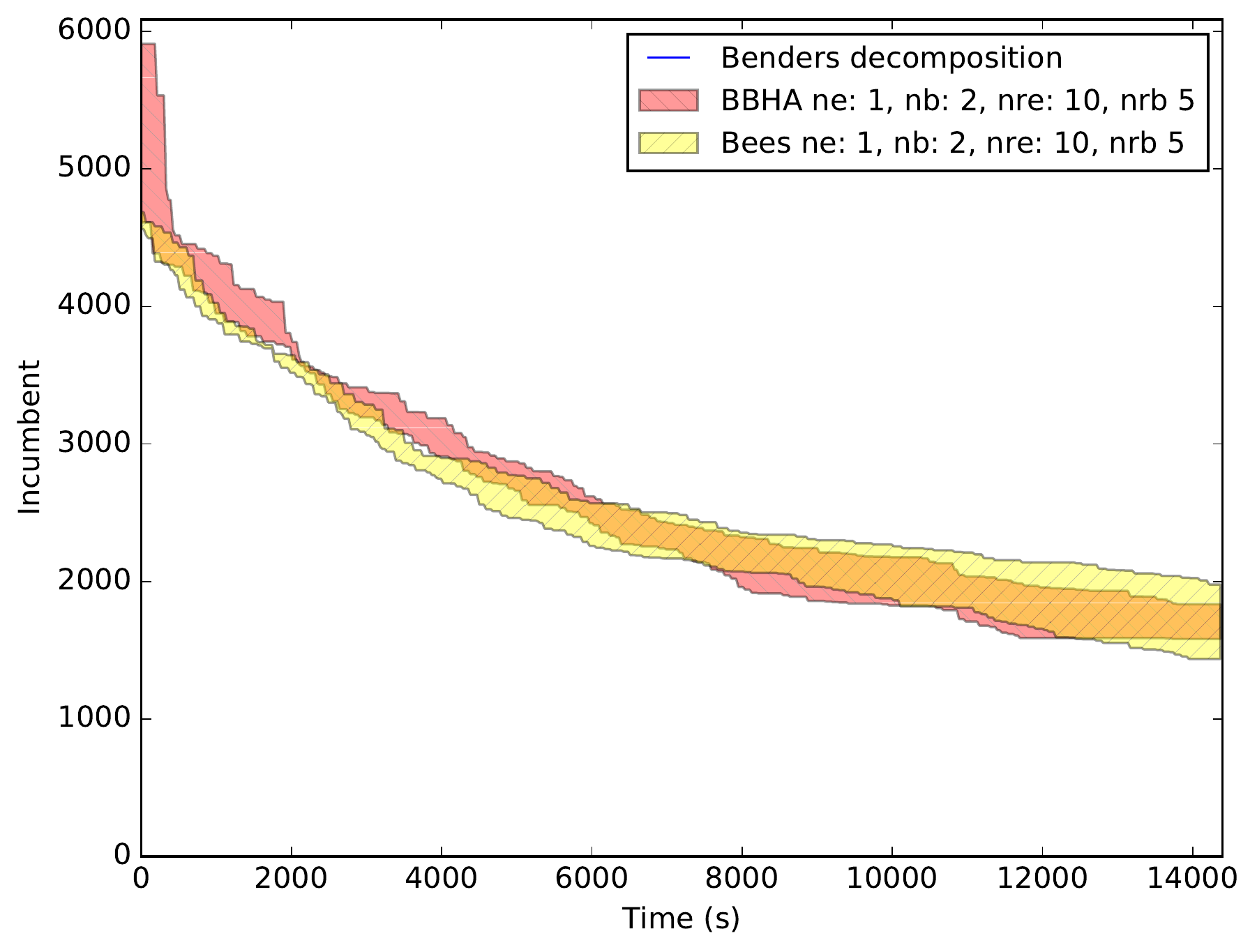}
  \caption{93-bus Long Peak scenario ne: 1, nb: 2, nre: 10, nrb: 5}
  \label{fig:93long_1_2_10_5}
\end{figure}

Where the problem favours the new approach, such as for the 93-bus
problem under the long peak scenario shown in Figure~\ref{fig:93long_1_2_10_5},
there is little discrepancy between the ranges of the BBHA and the BA although
the BA has better mean performance. Note: The Benders decomposition
incumbent value does not fall within the plotted range.

Like any other hybrid approach the BBHA is a compromise. A straight Benders
decomposition implementation running on the same computing infrastructure will
evaluate more of the search tree than the BBHA scout. Likewise, without the
continuously running scout or the trade off between producing Benders cuts and
local search the BA approach can dedicate more cores to evaluating candidate
solutions. As a result the straight Benders approach tends to outperform the
BBHA in terms of the lower bounds it produces. However since these are generally
of much less interest than finding good feasible solutions, and because the Bees
algorithm is unable to produce any lower bounds, we have restricted ourselves to
only compare upper bound solution quality in this analysis. 

However, empirically we have found the benefits of cut sharing largely negate
any compromise. In the first instance, the cuts generated by the Benders scout
improve the heuristic estimate of the fitness of the candidate solutions in the
worker solution pool. Likewise, the cuts generated in parallel by the
elite workers are typically in the neighbourhood of the incumbent solution
they prove useful to the Benders decomposition. Perhaps the clearest example of
this is shown in Figure~\ref{fig:46win_1_2_10_5}. Here, by sharing cuts between
workers and the simultaneous Benders decomposition, the BBHA is able to prove the optimal
solution faster than Benders decomposition alone, even though the Benders
decomposition has a resource advantage on the compute infrastructure. The effect
is also evident to a lesser extent in Figure~\ref{fig:46sum_1_2_10_5}.

\begin{figure}[!ht]
  \centering
  \includegraphics[width=.75\columnwidth]{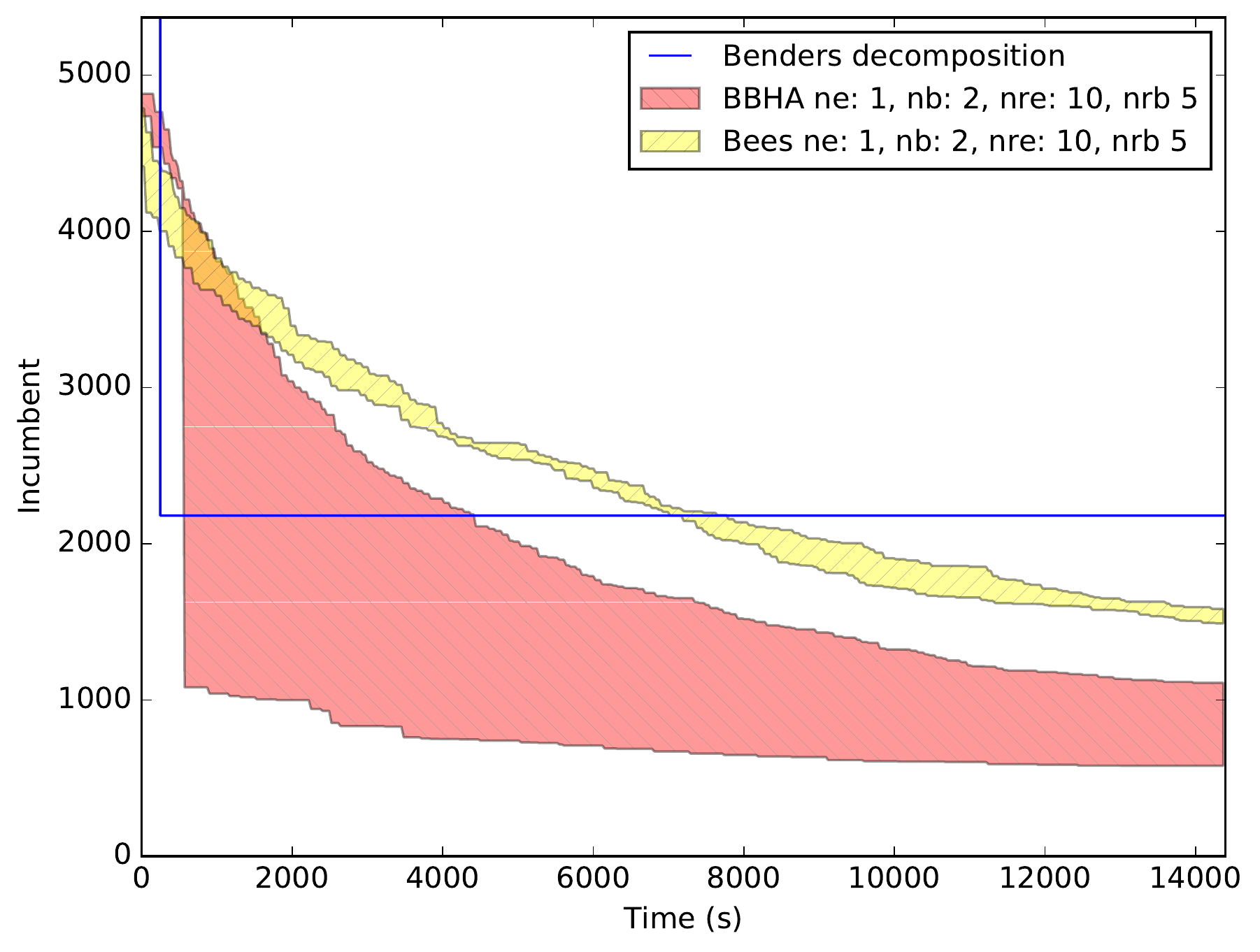}
  \caption{93-bus Short Peak scenario ne: 1, nb: 2, nre: 10, nrb: 5}
  \label{fig:93short_1_2_10_5}
\end{figure}

\begin{figure}[!ht]
  \centering
  \includegraphics[width=.75\columnwidth]{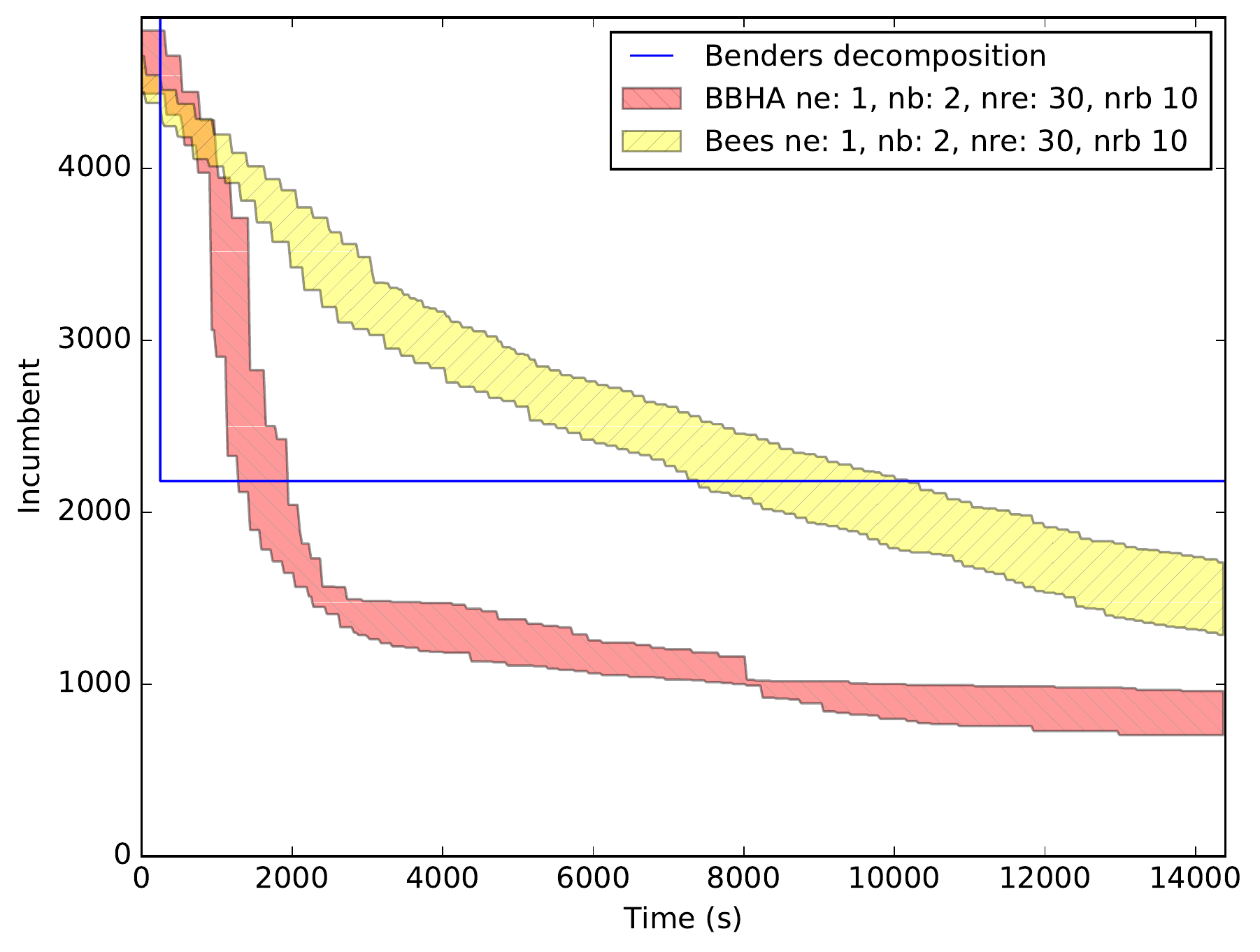}
  \caption{93-bus Short Peak scenario ne: 1, nb: 2, nre: 30, nrb: 10}
  \label{fig:93short_1_2_30_10}
\end{figure}

\begin{figure}[!ht]
  \centering
  \includegraphics[width=.75\columnwidth]{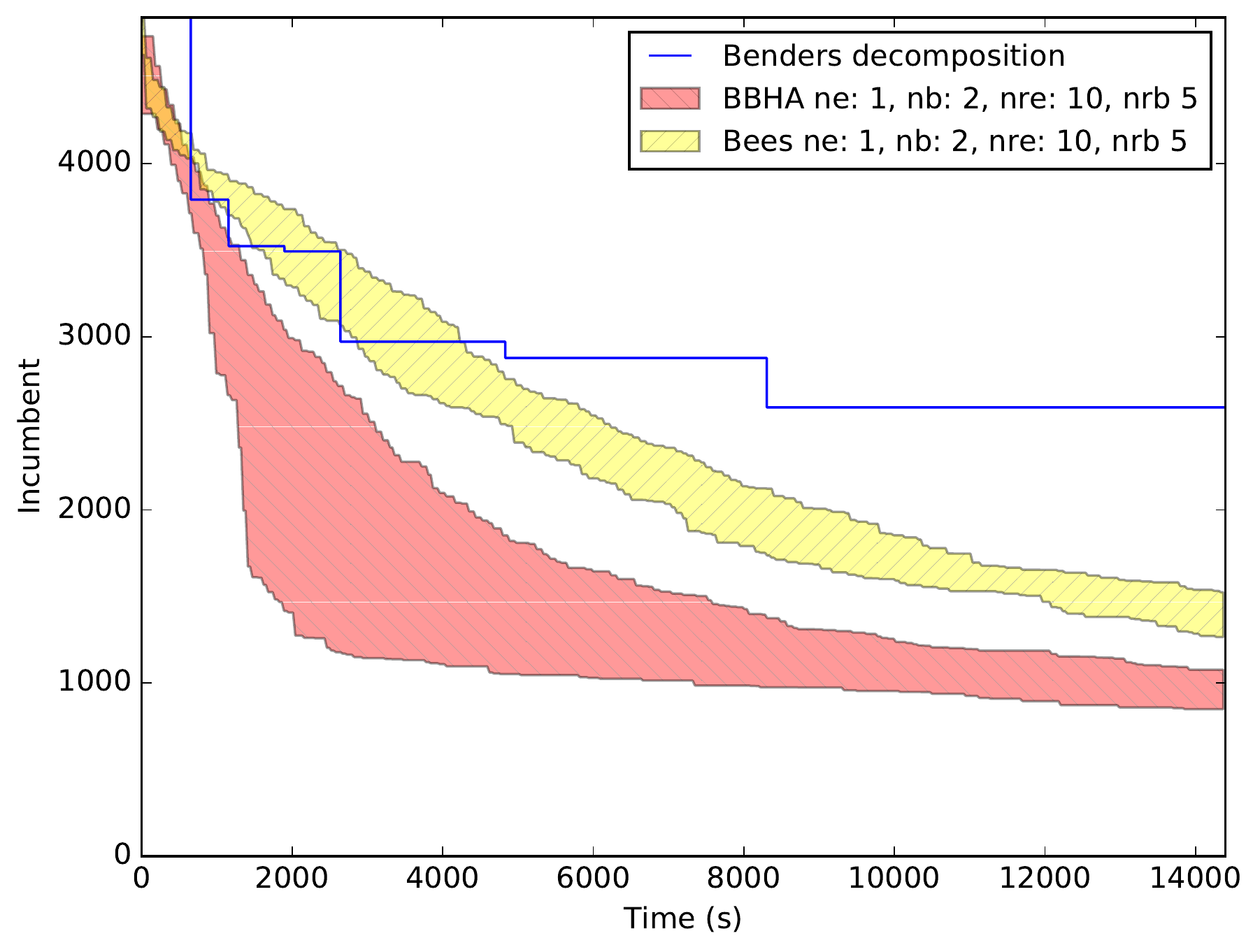}
  \caption{93-bus SGSC Summer scenario ne: 1, nb: 2, nre: 10, nrb: 5}
  \label{fig:93sum_1_2_10_5}
\end{figure}

Figure~\ref{fig:93short_1_2_10_5} shows the random effects of local search with
few workers, and the responsiveness of the BBHA to ``good'' Benders cuts.
These runs display a very large range of incumbent values because the worst of
the runs was unable to fully exploit the cut sharing. We observe similar random
effect in Figure~\ref{fig:93sum_1_2_10_5}. Here the BBHA shows sensitivity to
the parameters, and by increasing the workers available for local search
(parameter set: [ne: 1, nb: 2, nre: 30, nre: 10]) more of the search space is
evaluated each iteration and we observe less variance and broadly better.
Results across the 5 sample runs shown in Figure~\ref{fig:93short_1_2_30_10}.

Increasing the number of workers also significantly improves the
optimization of the 93-bus test system under the SGSC winter demand scenario
shown in Figures~\ref{fig:93win_1_2_10_5} (fewer workers) and
\ref{fig:93win_1_2_30_10} (more workers). In this case the entire range of
BBHA objective values improve upon that of the Benders decomposition by the end
of the optimization.


As noted in Section~\ref{subsec:cs_93-bus}, although the BBHA is a hybrid
matheuristic optimization technique, the use of Benders decomposition ensures
that the solution can be proven optimal if the algorithm is allowed to run for
sufficient time. If not allowed to run until completion, the library of cuts may
be used to produce a valid lower bound.

\section{Conclusion}
\label{sec:conclusion}

In this paper, we introduced a hybrid exact/meta-heuristic algorithm that
combined Benders decomposition and a Bees algorithm inspired approach. To the
best of our knowledge this is the first such matheuristic based on Benders
decomposition and the Bees algorithm.

The BBHA approach was demonstrated using a transmission network expansion and
energy storage planning model that is known to become more tractable when
decomposed into investment and operational subproblems.

The approach as been shown to combine the best performance of its component
parts in the segments of the problem domain where those parts excel, and to
improve upon the individual approaches where neither shows a substantial
advantage.

As the BBHA is general in nature and does not require any special problem
structure beyond the decomposition required by the Benders decomposition method.
The approach may be applied to any general decomposable mixed integer
programming problem. 

Future work will include a stochastic BBHA variant with multiple
probabilistically weighted subproblems. For problems such as the TESP
which rapidly become intractable as complexity and realism increases, a highly
parallelized stochastic algorithm is expected to significantly advance solution
quality.

\section*{Acknowledgment}
Melih Ozlen is supported by the Australian Research Council under the
Discovery Projects funding scheme (project~DP140104246).


\bibliographystyle{abbrvnat}
\bibliography{paper}

\begin{thebibliography}{52}
\providecommand{\natexlab}[1]{#1}
\providecommand{\url}[1]{\texttt{#1}}
\expandafter\ifx\csname urlstyle\endcsname\relax
  \providecommand{\doi}[1]{doi: #1}\else
  \providecommand{\doi}{doi: \begingroup \urlstyle{rm}\Url}\fi

\bibitem[Aghaei et~al.(2014)Aghaei, Muttaqi, Azizivahed, and
  Gitizadeh]{aghaei_distribution_2014}
J.~Aghaei, K.~M. Muttaqi, A.~Azizivahed, and M.~Gitizadeh.
\newblock Distribution expansion planning considering reliability and security
  of energy using modified {PSO} ({Particle} {Swarm} {Optimization}) algorithm.
\newblock \emph{Energy}, 65:\penalty0 398--411, Feb. 2014.
\newblock ISSN 0360-5442.
\newblock \doi{10.1016/j.energy.2013.10.082}.
\newblock URL
  \url{http://www.sciencedirect.com/science/article/pii/S0360544213009493}.
\newblock ADD.

\bibitem[Bahiense et~al.(2001)Bahiense, Oliveira, Pereira, and
  Granville]{bahiense_mixed_2001}
L.~Bahiense, G.~Oliveira, M.~Pereira, and S.~Granville.
\newblock A mixed integer disjunctive model for transmission network expansion.
\newblock \emph{IEEE Transactions on Power Systems}, 16\penalty0 (3):\penalty0
  560--565, 2001.
\newblock ISSN 0885-8950.
\newblock \doi{10.1109/59.932295}.

\bibitem[Benders(1962)]{benders_partitioning_1962}
J.~F. Benders.
\newblock Partitioning procedures for solving mixed-variables programming
  problems.
\newblock \emph{Numerische mathematik}, 4\penalty0 (1):\penalty0 238--252,
  1962.
\newblock URL \url{http://www.springerlink.com/index/g203830n1gm58w73.pdf}.

\bibitem[Berry et~al.(2013)Berry, Moore, Ward, Lindsay, and
  Proctor]{berry_national_2013}
A.~Berry, T.~Moore, J.~Ward, S.~Lindsay, and K.~Proctor.
\newblock National {Feeder} {Taxonomy} -- {Describing} a {Representative}
  {Feeder} {Set} for {Australian} {Electricity} {Distribution} {Networks}.
\newblock Technical report, CSIRO, Australia, June 2013.

\bibitem[Binato et~al.(2001)Binato, Pereira, and Granville]{binato_new_2001}
S.~Binato, M.~V.~F. Pereira, and S.~Granville.
\newblock A new {Benders} decomposition approach to solve power transmission
  network design problems.
\newblock \emph{IEEE Transactions on Power Systems}, 16\penalty0 (2):\penalty0
  235--240, 2001.
\newblock ISSN 0885-8950.
\newblock \doi{10.1109/59.918292}.

\bibitem[Boschetti et~al.(2009)Boschetti, Maniezzo, Roffilli, and
  R{\"o}hler]{boschetti_matheuristics:_2009}
M.~A. Boschetti, V.~Maniezzo, M.~Roffilli, and A.~B. R{\"o}hler.
\newblock Matheuristics: {Optimization}, {Simulation} and {Control}.
\newblock In \emph{Hybrid {Metaheuristics}}, Lecture {Notes} in {Computer}
  {Science}, pages 171--177. Springer, Berlin, Heidelberg, Oct. 2009.
\newblock ISBN 978-3-642-04917-0 978-3-642-04918-7.
\newblock \doi{10.1007/978-3-642-04918-7_13}.
\newblock URL
  \url{https://link-springer-com.ezproxy.lib.monash.edu.au/chapter/10.1007/978-3-642-04918-7_13}.

\bibitem[Chang(2016)]{chang_south_2016}
C.~Chang.
\newblock South {Australia} blackout highlights energy security.
\newblock \emph{News Limited}, Oct. 2016.
\newblock URL
  \url{http://www.news.com.au/technology/environment/can-we-rely-on-renewables/news-story/e727416fcba0bd866ebccb069af6255e}.

\bibitem[Chu and Beasley(1997)]{chu_genetic_1997}
P.~C. Chu and J.~E. Beasley.
\newblock A genetic algorithm for the generalised assignment problem.
\newblock \emph{Computers \& Operations Research}, 24\penalty0 (1):\penalty0
  17--23, Jan. 1997.
\newblock ISSN 0305-0548.
\newblock \doi{10.1016/S0305-0548(96)00032-9}.
\newblock URL
  \url{http://www.sciencedirect.com/science/article/pii/S0305054896000329}.

\bibitem[Clack et~al.(2015)Clack, Xie, and MacDonald]{clack_linear_2015}
C.~T.~M. Clack, Y.~Xie, and A.~E. MacDonald.
\newblock Linear programming techniques for developing an optimal electrical
  system including high-voltage direct-current transmission and storage.
\newblock \emph{International Journal of Electrical Power \& Energy Systems},
  68:\penalty0 103--114, June 2015.
\newblock ISSN 0142-0615.
\newblock \doi{10.1016/j.ijepes.2014.12.049}.
\newblock URL
  \url{http://www.sciencedirect.com/science/article/pii/S0142061514007765}.
\newblock r.

\bibitem[Costa(2005)]{costa_survey_2005}
A.~M. Costa.
\newblock A survey on benders decomposition applied to fixed-charge network
  design problems.
\newblock \emph{Computers \& Operations Research}, 32\penalty0 (6):\penalty0
  1429--1450, June 2005.
\newblock ISSN 0305-0548.
\newblock \doi{10.1016/j.cor.2003.11.012}.
\newblock URL
  \url{http://www.sciencedirect.com/science/article/pii/S0305054803003435}.

\bibitem[Crossland et~al.(2014)Crossland, Jones, and
  Wade]{crossland_planning_2014}
A.~F. Crossland, D.~Jones, and N.~S. Wade.
\newblock Planning the location and rating of distributed energy storage in
  {LV} networks using a genetic algorithm with simulated annealing.
\newblock \emph{International Journal of Electrical Power \& Energy Systems},
  59:\penalty0 103--110, July 2014.
\newblock ISSN 0142-0615.
\newblock \doi{10.1016/j.ijepes.2014.02.001}.
\newblock URL
  \url{http://www.sciencedirect.com/science/article/pii/S0142061514000532}.
\newblock r.

\bibitem[De~J~Silva et~al.(2005)De~J~Silva, Rider, Romero, Garcia, and
  Murari]{de_j_silva_transmission_2005}
I.~De~J~Silva, M.~Rider, R.~Romero, A.~Garcia, and C.~Murari.
\newblock Transmission network expansion planning with security constraints.
\newblock \emph{Generation, Transmission and Distribution, IEE Proceedings-},
  152\penalty0 (6):\penalty0 828--836, 2005.
\newblock ISSN 1350-2360.
\newblock \doi{10.1049/ip-gtd:20045217}.

\bibitem[Dilwali et~al.(2016)Dilwali, Gunnaasankaraan, Viswanath, and
  Mahata]{dilwali_transmission_2016}
K.~Dilwali, H.~Gunnaasankaraan, A.~Viswanath, and K.~Mahata.
\newblock Transmission expansion planning using benders decomposition and local
  branching.
\newblock In \emph{2016 {IEEE} {Power} and {Energy} {Conference} at {Illinois}
  ({PECI})}, pages 1--8, Feb. 2016.
\newblock \doi{10.1109/PECI.2016.7459265}.

\bibitem[Ekwue and Cory(1984)]{ekwue_transmission_1984}
A.~Ekwue and B.~Cory.
\newblock Transmission {System} {Expansion} {Planning} by {Interactive}
  {Methods}.
\newblock \emph{IEEE Transactions on Power Apparatus and Systems},
  PAS-103\penalty0 (7):\penalty0 1583--1591, July 1984.
\newblock ISSN 0018-9510.
\newblock \doi{10.1109/TPAS.1984.318637}.

\bibitem[Ernst(2010)]{ernst_hybrid_2010}
A.~T. Ernst.
\newblock A hybrid {Lagrangian} {Particle} {Swarm} {Optimization} {Algorithm}
  for the degree-constrained minimum spanning tree problem.
\newblock In \emph{{IEEE} {Congress} on {Evolutionary} {Computation}}, pages
  1--8, 2010.
\newblock \doi{http://dx.doi.org/10.1109/CEC.2010.5585939}.
\newblock 00009.

\bibitem[Escobar et~al.(2004)Escobar, Gallego, and
  Romero]{escobar_multistage_2004}
A.~H. Escobar, R.~A. Gallego, and R.~Romero.
\newblock Multistage and coordinated planning of the expansion of transmission
  systems.
\newblock \emph{IEEE Transactions on Power Systems}, 19\penalty0 (2):\penalty0
  735--744, May 2004.
\newblock ISSN 0885-8950.
\newblock \doi{10.1109/TPWRS.2004.825920}.

\bibitem[Eshragh et~al.(2011)Eshragh, Filar, and
  Nazari]{eshragh_projection-adapted_2011}
A.~Eshragh, J.~Filar, and A.~Nazari.
\newblock A {Projection}-{Adapted} {Cross} {Entropy} ({PACE}) method for
  transmission network planning.
\newblock \emph{Energy Systems}, 2\penalty0 (2):\penalty0 189--208, May 2011.
\newblock ISSN 1868-3967, 1868-3975.
\newblock \doi{10.1007/s12667-011-0033-x}.
\newblock URL \url{http://link.springer.com/article/10.1007/s12667-011-0033-x}.

\bibitem[Geoffrion(1972)]{geoffrion_generalized_1972}
A.~M. Geoffrion.
\newblock Generalized {Benders} decomposition.
\newblock \emph{Journal of Optimization Theory and Applications}, 10\penalty0
  (4):\penalty0 237--260, Oct. 1972.
\newblock ISSN 0022-3239, 1573-2878.
\newblock \doi{10.1007/BF00934810}.
\newblock URL \url{http://link.springer.com/article/10.1007/BF00934810}.

\bibitem[Gomez-Iglesias et~al.(2013)Gomez-Iglesias, Ernst, and
  Singh]{gomez-iglesias_scalable_2013}
A.~Gomez-Iglesias, A.~Ernst, and G.~Singh.
\newblock Scalable {Multi} {Swarm}-{Based} {Algorithms} with {Lagrangian}
  {Relaxation} for {Constrained} {Problems}.
\newblock In \emph{2013 12th {IEEE} {International} {Conference} on {Trust},
  {Security} and {Privacy} in {Computing} and {Communications} ({TrustCom})},
  pages 1073--1080, July 2013.
\newblock \doi{10.1109/TrustCom.2013.241}.
\newblock 00000.

\bibitem[Gu et~al.(2012)Gu, McCalley, and Ni]{gu_coordinating_2012}
Y.~Gu, J.~McCalley, and M.~Ni.
\newblock Coordinating {Large}-{Scale} {Wind} {Integration} and {Transmission}
  {Planning}.
\newblock \emph{IEEE Transactions on Sustainable Energy}, 3\penalty0
  (4):\penalty0 652--659, 2012.
\newblock ISSN 1949-3029.
\newblock \doi{10.1109/TSTE.2012.2204069}.

\bibitem[Haffner et~al.(2000)Haffner, Monticelli, Garcia, Mantovani, and
  Romero]{haffner_branch_2000}
S.~Haffner, A.~Monticelli, A.~Garcia, J.~Mantovani, and R.~Romero.
\newblock Branch and bound algorithm for transmission system expansion planning
  using a transportation model.
\newblock \emph{Generation, Transmission and Distribution, IEE Proceedings-},
  147\penalty0 (3):\penalty0 149--156, 2000.
\newblock ISSN 1350-2360.
\newblock \doi{10.1049/ip-gtd:20000337}.

\bibitem[Hu et~al.(2012)Hu, Zhang, and Li]{hu_transmission_2012}
Z.~Hu, F.~Zhang, and B.~Li.
\newblock Transmission expansion planning considering the deployment of energy
  storage systems.
\newblock In \emph{2012 {IEEE} {Power} and {Energy} {Society} {General}
  {Meeting}}, pages 1--6, 2012.
\newblock \doi{10.1109/PESGM.2012.6344575}.

\bibitem[Karaboga et~al.(2014)Karaboga, Gorkemli, Ozturk, and
  Karaboga]{karaboga2014}
D.~Karaboga, B.~Gorkemli, C.~Ozturk, and N.~Karaboga.
\newblock A comprehensive survey: artificial bee colony (abc) algorithm and
  applications.
\newblock \emph{The Artificial Intelligence Review}, 42\penalty0 (1):\penalty0
  21--57, 06 2014.
\newblock URL
  \url{https://search-proquest-com.ezproxy.lib.monash.edu.au/docview/1524243122?accountid=12528}.
\newblock Copyright - Springer Science+Business Media Dordrecht 2014; Last
  updated - 2015-12-22.

\bibitem[Keck et~al.(2019)Keck, Lenzen, Vassallo, and Li]{keck2019impact}
F.~Keck, M.~Lenzen, A.~Vassallo, and M.~Li.
\newblock The impact of battery energy storage for renewable energy power grids
  in australia.
\newblock \emph{Energy}, 2019.

\bibitem[KNIGHT(1972)]{knight197228}
U.~KNIGHT.
\newblock {CHAPTER} 3 - {SOME} {FREQUENTLY} {USED} {ANALYTICAL} {TECHNIQUES}.
\newblock In U.~KNIGHT, editor, \emph{Power Systems Engineering and
  Mathematics}, International Series of Monographs in Electrical Engineering,
  pages 28 -- 51. Pergamon, 1972.
\newblock ISBN 978-0-08-016603-2.
\newblock \doi{http://dx.doi.org/10.1016/B978-0-08-016603-2.50007-4}.
\newblock URL
  \url{http://www.sciencedirect.com/science/article/pii/B9780080166032500074}.

\bibitem[Kostikas and Fragakis(2004)]{kostikas_genetic_2004}
K.~Kostikas and C.~Fragakis.
\newblock Genetic {Programming} {Applied} to {Mixed} {Integer} {Programming}.
\newblock In M.~Keijzer, U.-M. O’Reilly, S.~Lucas, E.~Costa, and T.~Soule,
  editors, \emph{Genetic {Programming}}, number 3003 in Lecture {Notes} in
  {Computer} {Science}, pages 113--124. Springer Berlin Heidelberg, Apr. 2004.
\newblock ISBN 978-3-540-21346-8 978-3-540-24650-3.
\newblock \doi{10.1007/978-3-540-24650-3_11}.
\newblock URL
  \url{http://link.springer.com/chapter/10.1007/978-3-540-24650-3_11}.

\bibitem[Latorre et~al.(2003)Latorre, Cruz, Areiza, and
  Villegas]{latorre_classification_2003}
G.~Latorre, R.~Cruz, J.~Areiza, and A.~Villegas.
\newblock Classification of publications and models on transmission expansion
  planning.
\newblock \emph{IEEE Transactions on Power Systems}, 18\penalty0 (2):\penalty0
  938--946, 2003.
\newblock ISSN 0885-8950.
\newblock \doi{10.1109/TPWRS.2003.811168}.

\bibitem[Ma et~al.(1997)Ma, Shahidehpour, and Marwali]{ma_transmission_1997}
H.~Ma, S.~Shahidehpour, and M.~Marwali.
\newblock Transmission constrained unit commitment based on {Benders}
  decomposition.
\newblock In \emph{American {Control} {Conference}, 1997. {Proceedings} of the
  1997}, volume~4, pages 2263--2267 vol.4, June 1997.
\newblock \doi{10.1109/ACC.1997.608991}.

\bibitem[MacRae et~al.(2016)MacRae, Ernst, and Ozlen]{macrae_benders_2016}
C.~A.~G. MacRae, A.~T. Ernst, and M.~Ozlen.
\newblock A {Benders} decomposition approach to transmission expansion planning
  considering energy storage.
\newblock \emph{Energy}, 112:\penalty0 795 -- 803, 2016.
\newblock ISSN 0360-5442.
\newblock \doi{http://dx.doi.org/10.1016/j.energy.2016.06.080}.
\newblock URL
  \url{http://www.sciencedirect.com/science/article/pii/S0360544216308544}.

\bibitem[Marino et~al.(2018)Marino, Marufuzzaman, Hu, and Sarder]{Marino2018}
C.~Marino, M.~Marufuzzaman, M.~Hu, and M.~Sarder.
\newblock Developing a cchp-microgrid operation decision model under
  uncertainty.
\newblock \emph{Computers \& Industrial Engineering}, 115\penalty0 (Supplement
  C):\penalty0 354 -- 367, 2018.
\newblock ISSN 0360-8352.
\newblock \doi{https://doi.org/10.1016/j.cie.2017.11.021}.
\newblock URL
  \url{http://www.sciencedirect.com/science/article/pii/S0360835217305533}.

\bibitem[Moradi et~al.(2011)Moradi, Razi, and Fatahi]{moradi_application_2011}
S.~Moradi, P.~Razi, and L.~Fatahi.
\newblock On the application of bees algorithm to the problem of crack
  detection of beam-type structures.
\newblock \emph{Computers \& Structures}, 89\penalty0 (23--24):\penalty0
  2169--2175, Dec. 2011.
\newblock ISSN 0045-7949.
\newblock \doi{10.1016/j.compstruc.2011.08.020}.
\newblock URL
  \url{http://www.sciencedirect.com/science/article/pii/S0045794911002434}.

\bibitem[\"{O}zbakir et~al.(2010)\"{O}zbakir, Baykaso\u{g}lu, and
  Tapkan]{ozbakir_bees_2010}
L.~\"{O}zbakir, A.~Baykaso\u{g}lu, and P.~Tapkan.
\newblock Bees algorithm for generalized assignment problem.
\newblock \emph{Applied Mathematics and Computation}, 215\penalty0
  (11):\penalty0 3782--3795, Feb. 2010.
\newblock ISSN 0096-3003.
\newblock \doi{10.1016/j.amc.2009.11.018}.
\newblock URL
  \url{http://www.sciencedirect.com/science/article/pii/S0096300309010078}.

\bibitem[Pereira and Pinto(1991)]{Pereira1991}
M.~V.~F. Pereira and L.~M. V.~G. Pinto.
\newblock Multi-stage stochastic optimization applied to energy planning.
\newblock \emph{Mathematical Programming}, 52\penalty0 (1):\penalty0 359--375,
  May 1991.
\newblock ISSN 1436-4646.
\newblock \doi{10.1007/BF01582895}.
\newblock URL \url{https://doi.org/10.1007/BF01582895}.

\bibitem[Pereira et~al.(1985)Pereira, Pinto, Cunha, and
  Oliveira]{pereira_decomposition_1985}
M.~V.~F. Pereira, L.~M. V.~G. Pinto, S.~H.~F. Cunha, and G.~Oliveira.
\newblock A {Decomposition} {Approach} {To} {Automated}
  {Generation}/{Transmission} {Expansion} {Planning}.
\newblock \emph{IEEE Transactions on Power Apparatus and Systems},
  PAS-104\penalty0 (11):\penalty0 3074--3083, 1985.
\newblock ISSN 0018-9510.
\newblock \doi{10.1109/TPAS.1985.318815}.

\bibitem[Pham et~al.(2007)Pham, Koc, Lee, and Phrueksanant]{pham2007using}
D.~Pham, E.~Koc, J.~Lee, and J.~Phrueksanant.
\newblock Using the bees algorithm to schedule jobs for a machine.
\newblock In \emph{Proc eighth international conference on laser metrology, CMM
  and machine tool performance, LAMDAMAP, Euspen, UK, Cardiff}, pages 430--439,
  2007.

\bibitem[Pham and Castellani(2009)]{pham_bees_2009}
D.~T. Pham and M.~Castellani.
\newblock The {Bees} {Algorithm}: {Modelling} foraging behaviour to solve
  continuous optimization problems.
\newblock \emph{Proceedings of the Institution of Mechanical Engineers, Part C:
  Journal of Mechanical Engineering Science}, 223\penalty0 (12):\penalty0
  2919--2938, Dec. 2009.
\newblock ISSN 0954-4062, 2041-2983.
\newblock \doi{10.1243/09544062JMES1494}.
\newblock URL \url{http://pic.sagepub.com/content/223/12/2919}.

\bibitem[Poojari and Beasley(2009)]{poojari_improving_2009}
C.~A. Poojari and J.~E. Beasley.
\newblock Improving benders decomposition using a genetic algorithm.
\newblock \emph{European Journal of Operational Research}, 199\penalty0
  (1):\penalty0 89--97, Nov. 2009.
\newblock ISSN 0377-2217.
\newblock \doi{10.1016/j.ejor.2008.10.033}.
\newblock URL
  \url{http://www.sciencedirect.com/science/article/pii/S0377221708009740}.

\bibitem[Puchinger et~al.(2010)Puchinger, Raidl, and
  Pirkwieser]{puchinger_metaboosting2010}
J.~Puchinger, G.~R. Raidl, and S.~Pirkwieser.
\newblock {MetaBoosting}: {Enhancing} {Integer} {Programming} {Techniques} by
  {Metaheuristics}.
\newblock In V.~Maniezzo, T.~Stützle, and S.~Voß, editors,
  \emph{Matheuristics: {Hybridizing} {Metaheuristics} and {Mathematical}
  {Programming}}, Annals of {Information} {Systems}, pages 71--102. Springer
  US, Boston, MA, 2010.
\newblock ISBN 978-1-4419-1306-7.
\newblock \doi{10.1007/978-1-4419-1306-7_3}.
\newblock URL \url{https://doi.org/10.1007/978-1-4419-1306-7_3}.

\bibitem[Qiu et~al.(2017)Qiu, Xu, Wang, Dvorkin, and Kirschen]{Qiu2017}
T.~Qiu, B.~Xu, Y.~Wang, Y.~Dvorkin, and D.~Kirschen.
\newblock Stochastic multistage coplanning of transmission expansion and energy
  storage.
\newblock \emph{IEEE Transactions on Power Systems}, 32\penalty0 (1):\penalty0
  643--651, 2017.
\newblock \doi{10.1109/TPWRS.2016.2553678}.
\newblock URL
  \url{https://www.scopus.com/inward/record.uri?eid=2-s2.0-85008499528&doi=10.1109%2fTPWRS.2016.2553678&partnerID=40&md5=f166b17d283fd7430955494cb90ffd9a}.

\bibitem[Rebennack(2016)]{Rebennack2016}
S.~Rebennack.
\newblock Combining sampling-based and scenario-based nested benders
  decomposition methods: application to stochastic dual dynamic programming.
\newblock \emph{Mathematical Programming}, 156\penalty0 (1):\penalty0 343--389,
  Mar 2016.
\newblock ISSN 1436-4646.
\newblock \doi{10.1007/s10107-015-0884-3}.
\newblock URL \url{https://doi.org/10.1007/s10107-015-0884-3}.

\bibitem[Romero et~al.(2002)Romero, Monticelli, Garcia, and
  Haffner]{romero_test_2002}
R.~Romero, A.~Monticelli, A.~Garcia, and S.~Haffner.
\newblock Test systems and mathematical models for transmission network
  expansion planning.
\newblock \emph{Generation, Transmission and Distribution, IEE Proceedings-},
  149\penalty0 (1):\penalty0 27--36, 2002.
\newblock ISSN 1350-2360.
\newblock \doi{10.1049/ip-gtd:20020026}.

\bibitem[Saharidis et~al.(2010)Saharidis, Minoux, and
  Ierapetritou]{saharidis_accelerating_2010}
G.~K.~D. Saharidis, M.~Minoux, and M.~G. Ierapetritou.
\newblock Accelerating {Benders} method using covering cut bundle generation.
\newblock \emph{International Transactions in Operational Research},
  17\penalty0 (2):\penalty0 221--237, Mar. 2010.
\newblock ISSN 1475-3995.
\newblock \doi{10.1111/j.1475-3995.2009.00706.x}.
\newblock URL
  \url{http://onlinelibrary.wiley.com.ezproxy.lib.rmit.edu.au/doi/10.1111/j.1475-3995.2009.00706.x/abstract}.

\bibitem[Sedghi et~al.(2013)Sedghi, Aliakbar-Golkar, and
  Haghifam]{sedghi_distribution_2013}
M.~Sedghi, M.~Aliakbar-Golkar, and M.-R. Haghifam.
\newblock Distribution network expansion considering distributed generation and
  storage units using modified {PSO} algorithm.
\newblock \emph{International Journal of Electrical Power \& Energy Systems},
  52:\penalty0 221--230, Nov. 2013.
\newblock ISSN 0142-0615.
\newblock \doi{10.1016/j.ijepes.2013.03.041}.
\newblock URL
  \url{http://www.sciencedirect.com/science/article/pii/S0142061513001580}.
\newblock r*.

\bibitem[Sheikh et~al.(2015)Sheikh, Komaki, and Malakooti]{Sheikh2015}
S.~Sheikh, M.~Komaki, and B.~Malakooti.
\newblock Integrated risk and multi-objective optimization of energy systems.
\newblock \emph{Computers \& Industrial Engineering}, 90\penalty0 (Supplement
  C):\penalty0 1 -- 11, 2015.
\newblock ISSN 0360-8352.
\newblock \doi{https://doi.org/10.1016/j.cie.2015.08.008}.
\newblock URL
  \url{http://www.sciencedirect.com/science/article/pii/S0360835215003447}.

\bibitem[Sirikum et~al.(2007)Sirikum, Techanitisawad, and
  Kachitvichyanukul]{sirikum_new_2007}
J.~Sirikum, A.~Techanitisawad, and V.~Kachitvichyanukul.
\newblock A {New} {Efficient} {GA}-{Benders}' {Decomposition} {Method}: {For}
  {Power} {Generation} {Expansion} {Planning} {With} {Emission} {Controls}.
\newblock \emph{IEEE Transactions on Power Systems}, 22\penalty0 (3):\penalty0
  1092--1100, Aug. 2007.
\newblock ISSN 0885-8950.
\newblock \doi{10.1109/TPWRS.2007.901092}.

\bibitem[Sorokin et~al.(2012)Sorokin, Portela, and
  Pardalos]{sorokin_algorithms_2012}
A.~Sorokin, J.~Portela, and P.~Pardalos.
\newblock Algorithms and {Models} for {Transmission} {Expansion} {Planning}.
\newblock In A.~Sorokin, S.~Rebennack, P.~M. Pardalos, N.~A. Iliadis, and
  M.~V.~F. Pereira, editors, \emph{Handbook of {Networks} in {Power} {Systems}
  {I}}, Energy {Systems}, pages 395--433. Springer Berlin Heidelberg, 2012.
\newblock ISBN 978-3-642-23192-6.

\bibitem[Sum-Im et~al.(2009)Sum-Im, Taylor, Irving, and
  Song]{sum-im_differential_2009}
T.~Sum-Im, G.~A. Taylor, M.~R. Irving, and Y.~H. Song.
\newblock Differential evolution algorithm for static and multistage
  transmission expansion planning.
\newblock \emph{IEE proceedings. Generation, transmission, and distribution.},
  3\penalty0 (4):\penalty0 365--384, 2009.
\newblock ISSN 13502360.
\newblock URL
  \url{http://search.proquest.com.ezproxy.lib.rmit.edu.au/docview/1627082886/abstract/CDCD7BBBBC5145A6PQ/1}.

\bibitem[Thiruvady et~al.(2013)Thiruvady, Singh, Ernst, and
  Meyer]{thiruvady_constraint-based_2013}
D.~Thiruvady, G.~Singh, A.~T. Ernst, and B.~Meyer.
\newblock Constraint-based {ACO} for a shared resource constrained scheduling
  problem.
\newblock \emph{International Journal of Production Economics}, 141\penalty0
  (1):\penalty0 230--242, 2013.

\bibitem[Thiruvady et~al.(2014)Thiruvady, Ernst, and
  Wallace]{thiruvady_lagrangian-aco_2014}
D.~Thiruvady, A.~Ernst, and M.~Wallace.
\newblock A {Lagrangian}-{ACO} matheuristic for car sequencing.
\newblock \emph{EURO Journal on Computational Optimization}, 2\penalty0
  (4):\penalty0 279--296, 2014.

\bibitem[Wood(2012)]{wood_integrating_2012}
J.~Wood.
\newblock Integrating renewables into the grid: {Applying} {UltraBattery}
  \#x00ae; {Technology} in {MW} scale energy storage solutions for continuous
  variability management.
\newblock In \emph{2012 {IEEE} {International} {Conference} on {Power} {System}
  {Technology} ({POWERCON})}, pages 1--4, Oct. 2012.
\newblock \doi{10.1109/PowerCon.2012.6401258}.

\bibitem[Wu and Shahidehpour(2010)]{wu_accelerating_2010}
L.~Wu and M.~Shahidehpour.
\newblock Accelerating the {Benders} decomposition for network-constrained unit
  commitment problems.
\newblock \emph{Energy Systems}, 1\penalty0 (3):\penalty0 339--376, Aug. 2010.
\newblock ISSN 1868-3975.
\newblock \doi{10.1007/s12667-010-0015-4}.
\newblock URL \url{https://doi.org/10.1007/s12667-010-0015-4}.

\bibitem[Zhao et~al.(2015)Zhao, Wu, Hu, Xu, and Rasmussen]{zhao_review_2015}
H.~Zhao, Q.~Wu, S.~Hu, H.~Xu, and C.~N. Rasmussen.
\newblock Review of energy storage system for wind power integration support.
\newblock \emph{Applied Energy}, 137:\penalty0 545--553, Jan. 2015.
\newblock ISSN 0306-2619.
\newblock \doi{10.1016/j.apenergy.2014.04.103}.
\newblock URL
  \url{http://www.sciencedirect.com/science/article/pii/S0306261914004668}.

\end{thebibliography}

\begin{table}[!h]
  \renewcommand{\arraystretch}{1.3}
  \fontsize{8}{10}\selectfont
  \caption{Parameters for the BBHA.}
  \label{tab:params}
  \centering
  \begin{tabular}{@{} c>{\raggedright}p{20em}c@{}}
    \toprule
    Name & Description & Default value \\
    \midrule
    ne & number of elite sites & 1 \\
    nb & number of best sites & 2 \\
    nre & recruited bees for elite sites & 10 \\
    nrb & recruited bees for remaining best sites & 5\\
    ngh & maximum size of neighbourhood for local search & 8 \\
    \bottomrule
  \end{tabular}
\end{table}

\begin{table}
  \centering
  \fontsize{8}{10}\selectfont
\begin{filecontents*}{tuningzz.txt}
obj,total_time,load,trapz_bees,ts,nre,scaled_trapz,nrb,ne,nb,iterations
44702.56,1808.116,Long peak,56426168.1292,IEEE-25-bus,10,0.260706648048,5,2,3,93
44702.56,1817.123,Long peak,65617611.9922,IEEE-25-bus,10,0.303174010261,5,1,4,92
43812.33,1817.285,Long peak,68383878.168,IEEE-25-bus,30,0.315955030242,10,1,2,56
44702.56,1820.305,Long peak,79659571.554,IEEE-25-bus,20,0.368052280943,15,1,4,39
44702.56,1814.477,Long peak,86350958.0494,IEEE-25-bus,10,0.398968591616,5,1,2,133
44702.56,1816.46,Long peak,88583191.7088,IEEE-25-bus,30,0.409282213368,25,1,3,32
55282.15,1811.793,Long peak,88588750.958,IEEE-25-bus,20,0.409307898848,15,1,2,67
43812.33,1853.469,Long peak,94731861.2494,IEEE-25-bus,30,0.437691000976,15,1,2,50
44702.56,1822.816,Long peak,97998429.4859,IEEE-25-bus,20,0.452783573869,10,1,3,57
45377.9,1820.465,Long peak,100429841.0,IEEE-25-bus,20,0.464017459968,10,3,4,36
44702.56,1905.349,Long peak,101580226.323,IEEE-25-bus,30,0.469332602065,20,1,3,36
45377.9,1805.288,Long peak,107690414.692,IEEE-25-bus,10,0.497563594553,5,1,3,88
44478.76,1849.929,Long peak,108654300.786,IEEE-25-bus,20,0.50201705154,5,2,3,48
43812.33,1823.841,Long peak,109004795.251,IEEE-25-bus,20,0.503636446233,10,2,3,52
46363.83,1819.196,Long peak,117213841.013,IEEE-25-bus,20,0.541564820161,5,3,4,30
56453.57,1883.395,Long peak,121898368.801,IEEE-25-bus,30,0.563208812265,10,3,4,22
55088.22,1827.963,Long peak,125544290.832,IEEE-25-bus,20,0.580054119029,15,1,3,53
55088.22,1811.475,Long peak,126006557.53,IEEE-25-bus,20,0.582189936598,10,1,2,63
44708.54,1842.026,Long peak,130231519.388,IEEE-25-bus,30,0.601710589526,20,2,3,30
44478.76,1835.736,Long peak,130857867.906,IEEE-25-bus,30,0.604604516725,15,1,3,44
51489.11,1845.344,Long peak,134264314.974,IEEE-25-bus,30,0.620343373828,20,3,4,24
43812.33,1800.211,Long peak,134593016.532,IEEE-25-bus,10,0.62186207843,5,3,4,69
44702.56,1822.615,Long peak,134996771.254,IEEE-25-bus,30,0.623727552263,20,1,2,50
44702.56,1898.531,Long peak,135723269.411,IEEE-25-bus,30,0.627084202304,15,2,3,32
55120.22,1833.296,Long peak,138238938.019,IEEE-25-bus,30,0.638707382682,15,3,4,24
48879.54,1831.608,Long peak,138485421.956,IEEE-25-bus,30,0.639846216017,10,1,4,40
51784.49,1837.557,Long peak,139784923.776,IEEE-25-bus,30,0.645850323238,15,1,4,35
43812.33,1814.051,Long peak,141840995.772,IEEE-25-bus,20,0.655350022684,15,3,4,30
46044.54,1821.33,Long peak,143090295.901,IEEE-25-bus,30,0.661122182304,10,1,3,47
55282.15,1838.311,Long peak,148156126.646,IEEE-25-bus,30,0.684527914021,10,2,3,35
60263.09,1812.03,Long peak,150194718.389,IEEE-25-bus,20,0.6939468492,15,2,3,34
50091.63,1801.222,Long peak,164907441.915,IEEE-25-bus,30,0.761924260415,25,1,4,22
45413.91,1811.332,Long peak,167706902.918,IEEE-25-bus,20,0.774858650942,10,1,4,38
52760.14,1850.153,Long peak,216435478.541,IEEE-25-bus,30,1.0,20,1,4,20
\end{filecontents*}
\pgfplotstabletypeset[
  columns={load, ne, nb, nre, nrb, obj, iterations, scaled_trapz},
  every head row/.style={
    before row={
      \toprule
  },
  after row={
    & & & & & (US\$10\textsuperscript{3}) & & \\
    \hline
  }
},
columns/load/.style={column name=Scenario, string type, column type=l},
columns/ne/.style={column name=ne},
columns/nb/.style={column name=nb},
columns/nre/.style={column name=nre},
columns/nrb/.style={column name=nrb},
columns/obj/.style={column name=Objective, precision=0, zerofill, std=-5,
  preproc/expr={##1/1}},
columns/interations/.style={column name=Iterations},
columns/scaled_trapz/.style={column name=Scaled Trapz, precision=2, zerofill, std=-5},
every last row/.style={
  after row=\bottomrule
},
  col sep=comma]
  {tuningzz.txt}
  \caption{Tuning results for IEEE 25-bus test system. Note: We use a composite
  trapezoidal rule to integrate along the time axis and then rescale (Scaled
  Trapz) which gives a measure of average solution quality over the half
  hour duration of the run. Lower is better.}
  \label{tab:tuning}
\end{table}

\begin{table}
  \centering
  \fontsize{7}{7}\selectfont
\begin{filecontents*}{tuning_paper_plots_09_sep.txt}
nat_mean,nat_worst,bee_mean,benders,bee_worst,bee_best,ps,ts,nat_best,ld
170644.483333,216458.67,107205.36,111840.23131,121394.82,100110.63,1 2 10 5,46-bus,144749.27,Long peak
128806.826667,142615.21,111251.96,111840.23131,127617.19,100110.63,1 2 30 10,46-bus,115754.15,Long peak
201617.23,241328.47,113597.683333,111840.23131,119453.01,110321.02,2 3 10 5,46-bus,180162.04,Long peak
113529.473333,138245.85,72355.41,72355.40797,72355.41,72355.41,1 2 10 5,46-bus,98384.7,Short peak
108905.216667,116960.58,72355.41,72355.40797,72355.41,72355.41,1 2 30 10,46-bus,102433.18,Short peak
145689.3,165592.33,72355.41,72355.40797,72355.41,72355.41,2 3 10 5,46-bus,125070.37,Short peak
65746.0166667,72356.25,46434.71,46434.7145,46434.71,46434.71,1 2 10 5,46-bus,52702.91,SGSC summer
55531.4133333,59616.67,46434.71,46434.7145,46434.71,46434.71,1 2 30 10,46-bus,48323.52,SGSC summer
93496.4433333,108743.6,46434.71,46434.7145,46434.71,46434.71,2 3 10 5,46-bus,78362.36,SGSC summer
94447.6666667,94794.21,59952.72,59952.72277,59952.72,59952.72,1 2 10 5,46-bus,93841.26,SGSC winter
77248.63,85922.72,59952.72,59952.72277,59952.72,59952.72,1 2 30 10,46-bus,68314.38,SGSC winter
101681.933333,118503.15,59952.72,59952.72277,59952.72,59952.72,2 3 10 5,46-bus,83985.08,SGSC winter
1705.92,1977.77,1743.23,9537.88793,1834.13,1581.22,1 2 10 5,93-bus,1429.31,Long peak
1723.29,1927.81,2120.86333333,9537.88793,2434.92,1891.62,1 2 30 10,93-bus,1550.24,Long peak
1832.58,1882.27,2258.85,9537.88793,2391.34,2137.23,2 3 10 5,93-bus,1740.81,Long peak
1526.77,1561.26,758.886666667,2181.05035,1110.53,579.23,1 2 10 5,93-bus,1490.24,Short peak
1539.22333333,1707.68,835.976666667,2181.05035,960.64,704.98,1 2 30 10,93-bus,1286.71,Short peak
1674.09333333,1792.89,1066.31,2181.05035,1425.75,829.03,2 3 10 5,93-bus,1582.94,Short peak
1403.26333333,1519.32,928.026666667,2592.20058,1076.99,848.93,1 2 10 5,93-bus,1266.19,SGSC summer
1421.82333333,1524.21,1039.27,2592.20058,1189.71,828.64,1 2 30 10,93-bus,1305.48,SGSC summer
1639.30333333,1733.61,1078.31333333,2592.20058,1093.63,1059.19,2 3 10 5,93-bus,1530.99,SGSC summer
1456.25666667,1507.05,1097.65666667,1077.04223,1444.57,897.01,1 2 10 5,93-bus,1360.55,SGSC winter
1558.55333333,1636.75,812.52,1077.04223,904.47,738.34,1 2 30 10,93-bus,1474.41,SGSC winter
1764.57,1853.96,1302.51333333,1077.04223,1501.19,1143.9,2 3 10 5,93-bus,1704.73,SGSC winter

\end{filecontents*}
\pgfplotstabletypeset[
  columns={ts, ld, ps, bee_worst, bee_mean, bee_best, nat_worst, nat_mean, nat_best, benders},
  every head row/.style={
    before row={
      \toprule
  },
  after row={
    & & & (US\$10\textsuperscript{3}) & (US\$10\textsuperscript{3}) & (US\$10\textsuperscript{3}) & (US\$10\textsuperscript{3}) & (US\$10\textsuperscript{3}) & (US\$10\textsuperscript{3}) & (US\$10\textsuperscript{3}) \\
    \hline
  }
},
columns/ts/.style={column name=Network, string type, column type=l},
columns/ld/.style={column name=Scenario, string type, column type=l},
columns/ps/.style={column name=Params, string type, column type=l},
columns/bee_worst/.style={column name=BBHA worst, fixed, precision=2, zerofill, column type = {r}},
columns/bee_mean/.style={column name=BBHA mean, fixed, precision=2, zerofill, column type = {r}},
columns/bee_best/.style={column name=BBHA best, fixed, precision=2, zerofill, column type = {r}},
columns/nat_worst/.style={column name=Bee worst, fixed, precision=2, zerofill, column type = {r}},
columns/nat_mean/.style={column name=Bee mean, fixed, precision=2, zerofill, column type = {r}},
columns/nat_best/.style={column name=Bee best, fixed, precision=2, zerofill, column type = {r}},
columns/benders/.style={column name=Benders, fixed, precision=2, zerofill, column type = {r}},
every last row/.style={
  after row=\bottomrule
},
  col sep=comma]
  {tuning_paper_plots_09_sep.txt}
  \caption{Results for the BBHA compared with the basic Bee Algorithm and
      Benders Decompostion on their own. The BBHA \& Bee algorithm
      were run three times with different parameter settings for each instance.} 
  \label{tab:results}
\end{table}

\end{document}